\documentclass[12pt]{amsart}
\usepackage{amsmath,amsthm,upref,amssymb}
\usepackage[dvips]{graphicx}
\usepackage{psfrag}
\usepackage{wrapfig}

\newcommand*{\newconcept}[1]{\textbf{#1}}

\newcommand*\N{\mathbf N}
\newcommand*\Z{\mathbf Z}
\newcommand*\Q{\mathbf Q}
\newcommand*\R{\mathbf R}
\newcommand*\C{\mathbf C}
\newcommand*\HH{\mathbf H}

\newcommand*\cH{{\mathcal{H}}}

\newcommand*\cF{{\mathcal{F}}}

\newcommand*\cU{{\mathcal{U}}}
\newcommand*\cC{{\mathcal{C}}}

\DeclareMathOperator{\Card}{Card}

\DeclareMathOperator{\vol}{vol}
\DeclareMathOperator\Aff{Aff}
\DeclareMathOperator\Diff{Diff}
\DeclareMathOperator\SL{SL}

\renewcommand*\Re{\operatorname{Re}}
\renewcommand*\Im{\operatorname{Im}}
\renewcommand*\d{\mathrm{d}}
\newcommand*\Id{\mathrm{Id}}
\newcommand*\U{U} 
\newcommand*\V{R} 
\newcommand*{\lquotient}[2]{\left.\raisebox{-0.2ex}{$#2$}%
  \backslash\raisebox{0.2ex}{$#1$}\right.}
\newcommand*{\rquotient}[2]{\left.\raisebox{0.2ex}{$#1$}%
  /\raisebox{-0.2ex}{$#2$}\right.}
\renewcommand*{\pmod}[1]{\;[#1]}

\newcommand\mat[4]{%
  \bigl( \begin{smallmatrix} #1&#3\\ #2&#4
  \end{smallmatrix}\bigr)
  }

\newtheorem{Theorem}{Theorem}[section]
\newtheorem*{Theorem*}{Theorem}
\newtheorem{Proposition}[Theorem]{Proposition}
\newtheorem{Lemma}[Theorem]{Lemma}
\newtheorem{Corollary}[Theorem]{Corollary}
\newtheorem{Conjecture}[Theorem]{Conjecture}
\newtheorem*{GBFormula}{The Gauss--Bonnet Formula}

\newtheorem*{Veechdichotomy}{Theorem (Veech dichotomy)}
\theoremstyle{remark}
\newtheorem*{Remark}{Remark}
\newtheorem*{Remarks}{Remarks}

\begin{document}

\title
{Prime arithmetic Teich\-m\"uller discs in $\cH(2)$}

%
\author{Pascal Hubert}
\address{%
IML, UMR CNRS 6206, Universit\'e de la M\'editerran\'ee,
Campus de Luminy, case 907, 13288 Marseille cedex 9, France.%
}
\email{hubert@iml.univ-mrs.fr}
\author{Samuel Leli\`evre}
\address{%
\textsc{Irmar}, UMR CNRS 6625, Universit\'e de Rennes~1,
Campus Beaulieu, 35042 Rennes cedex, France;%
\newline
\indent
I3M, UMR CNRS 5149, Universit\'e Montpellier~2, case 51,
Place Eug\`ene Bataillon, 34095 Montpellier cedex 5, France;%
\newline
\indent
IML, UMR CNRS 6206, Universit\'e de la M\'editerran\'ee,
Campus de Luminy, case 907, 13288 Marseille cedex 9, France.%
}
\email{samuel.lelievre@polytechnique.org}
\urladdr{http://carva.org/samuel.lelievre/}

\date{26 October 2004.}

\begin{abstract}
It is well-known that Teich\-m\"uller discs that pass through
``integer points'' of the moduli space of abelian differentials are
very special: they are closed complex geodesics.
However, the structure of these special Teich\-m\"uller discs is
mostly unexplored: their number, genus, area, cusps, etc.

We prove that in genus two all translation surfaces in $\cH(2)$ tiled
by a prime number $n > 3$ of squares fall into exactly two
Teich\-m\"uller discs, only one of them with elliptic points, and that
the genus of these discs has a cubic growth rate in $n$.
\end{abstract}

\maketitle
\thispagestyle{empty}



\setcounter{tocdepth}{1}
\tableofcontents
\normalsize

%

\section{Introduction}
\label{sec:intro}
In his fundamental paper of 1989, Veech studied the finite-volume
Teich\-m\"uller discs.
Translation surfaces with such discs, called Veech surfaces, enjoy
very interesting dynamical properties: their directional flows are
either completely periodic or uniquely ergodic.
An abundant literature exists on Veech surfaces: Veech
\cite{Ve89,Ve92}, Gutkin--Judge \cite{GuJu96,GuJu00}, Vorobets
\cite{Vo}, Ward \cite{Wa}, Kenyon--Smillie \cite{KeSm},
Hubert--Schmidt \cite{HuSc00,HuSc01}, Gutkin--Hubert--Schmidt
\cite{GuHuSc}, Calta \cite{Ca}, McMullen \cite{Mc}\ldots

The simplest examples of Veech surfaces are translation covers of the
torus (ramified over a single point), called square-tiled surfaces.
They are those translation surfaces whose stabilizer in $\SL(2,\R)$ is
arithmetic (commensurable with $\SL(2,\Z)$), by a theorem of Gutkin
and Judge.
These surfaces (and many more!)  were introduced by
Thurston \cite{Th} and studied on the dynamical aspect by Gutkin
\cite{Gu}, Veech \cite{Ve87} and Gutkin--Judge \cite{GuJu96,GuJu00}.
Square-tiled surfaces can be viewed as the ``integer points'' of the
moduli spaces of holomorphic $1$-forms.
The asymptotic number of integer points in a large ball was used by
Zorich \cite{Zo} and Eskin--Okounkov \cite{EsOk} to compute volumes of
strata of abelian differentials.

It was known for years that Teich\-m\"uller discs passing through
these integer points in the moduli space are very special: they are
closed (\emph{complex}) geodesics.
Despite enormous interest to invariant submanifolds (especially to the
simplest ones: those of complex dimension one), absolutely nothing was
known about the structure of these special Teich\-m\"uller discs:
about their number, genus, area, cusps, etc.
It was neither known which ``integer points'' belong to the same
Teich\-m\"uller disc.

\subsection{Main results}

In this paper, we study square-tiled surfaces in the stratum $\cH(2)$.
This stratum is the moduli space of holomorphic $1$-forms with a
unique (double) zero on a surface of genus two.
For surfaces tiled by a prime number of squares, we show:

\begin{Theorem}
\label{thm:twoorbits}
For any prime $n\geqslant 5$, the $\SL(2,\R)$-orbits of
$n$-square-tiled surfaces in $\cH(2)$ form two Teich\-m\"uller discs
$D_A(n)$ and $D_B(n)$.
\end{Theorem}

\begin{Theorem}
\label{thm:countings}
$D_A(n)$ and $D_B(n)$ can be seen as the unit tangent bundles to
orbifold surfaces with the following asymptotic behavior:\\
\textbullet\enspace genus $\sim c \, n^3$, with $c_A = c_B =
(3/16)(1/12)$,\\
\textbullet\enspace area $\sim c \, n^3$, with $c_A = c_B =
(3/16)(\pi/3)$,\\
\textbullet\enspace number of cusps $\sim c \, n^2$, with $c_A = 1/24$
and $c_B = 1/8$,\\
\textbullet\enspace number of elliptic points $O(n)$, one of them having none.
\end{Theorem}

\begin{Proposition}
\label{prn:L}
All these discs arise from L-shaped billiards.
\end{Proposition}

Our results are extended by McMullen \cite{Mc2} to describe the
repartition into different orbits of \emph{all} Veech surfaces in
$\cH(2)$.
In particular, the invariant introduced in \S\,\ref{sec:invariant}
also determines orbits in the nonprime case.

\subsection{Side results}

We find the following as side results of our study:

\smallskip

\textbullet\enspace \textbf{One-cylinder directions.}
\begin{Proposition}
\label{prn:one:cyl:dir}
All surfaces in $\cH(2)$ tiled by a prime number of squares have
one-cylinder directions i.e.\ directions in which they decompose into
one single cylinder.
\end{Proposition}

\smallskip

\textbullet\enspace \textbf{Discs without elliptic points.}
During some time, the search for new Veech surfaces focused on
examples arising from billiards in rational-angled polygons.
Angles of the billiard table not multiples of the right angle lead to
elliptic elements in the Veech group.
Billiards with all angles multiples of the right angle have however
recently been studied, especially L-shaped billiards (see \cite{Mc}).

\smallskip

\textbullet\enspace \textbf{Discs of (arbitrary high) positive genus.}
When a Veech group has positive genus, the subgroup generated by its
parabolic elements has infinite index, and cannot be a lattice.
This implies that the naive algorithm which consists in finding
parabolic elements in the Veech group cannot lead to obtain the
\emph{whole} group not even up to finite index.

The surfaces arising from billiards in the regular polygons, studied
by Veech in \cite{Ve92}, have genus tending to infinity, and one could
probably show that the genus of their Veech groups also tends to
infinity, though Veech does not state this explicitly.

Our examples give families of Teich\-m\"uller discs of arbitrarily high
genus, the translation surfaces in these discs staying in genus two.

\smallskip

\textbullet\enspace \textbf{Noncongruence subgroups.}
Since we deal with families of subgroups of $\SL(2,\Z)$, it is natural
to check whether they belong to the well-known family of congruence
subgroups.
Appendix~\ref{app:n:3:5} provides an example of a Veech group that is
a non-congruence subgroup of $\SL(2,\Z)$.
Another example was given by G.~Schmith\"usen \cite{Schmi}.
A detailed discussion of the congruence problem in this setting will
appear in \cite{HL}.

\smallskip

\textbullet\enspace \textbf{Deviation from the mean order.}
\begin{Proposition}
\label{prn:dev:mean:order}
The number of $n$-square-tiled surfaces in $\cH(2)$ for prime $n$ is
asymptotically $1/\zeta(4)$ times the mean order of the number of
$n$-square-tiled surfaces in $\cH(2)$.
\end{Proposition}

\subsection{Methods}

We parametrize square-tiled surfaces in $\cH(2)$ by using separatrix
diagrams as in \cite{KoZo}, \cite{Zo} and \cite{EsMaSc}.
These coordinates bring the study of Teich\-m\"uller discs of
$n$-square-tiled surfaces down to a combinatorial problem.

We want to describe the $\SL(2,\Z)$ orbits of these surfaces.
Using the fact that $\cH(2)$ is a hyperelliptic stratum, the
combinatorial representation of Weierstrass points allows us to show
there are at least two orbits for odd $n\geqslant 5$.
Showing there are only two is done for prime $n$ in a combinatorial
way, by a careful study of the action of generators of $\SL(2,\Z)$ on
square-tiled surfaces.

For the countings, we use generating functions.

\subsection{Related works}

Our counting results are very close to the formulae in \cite{EsMaSc}.
Eskin--Masur--Schmoll calculate Siegel--Veech constants for torus
coverings in genus two.
In $\cH(2)$, these calculations are based on counting the square-tiled
surfaces with a given number of squares.
The originality of our work is to count square-tiled surfaces disc by
disc.

There are also analogies with Schmoll's work \cite{Schmo}.
He computes the explicit Veech groups of tori with two marked points
and the quadratic asymptotics for theses surfaces.
Some of the methods he uses are intimately linked to those used in our
work.
The Veech groups he exhibits are all congruence subgroups.

A computer program allows to give all the geometric information on
Teich\-m\"uller discs of square-tiled surfaces in $\cH(2)$.
Schmith\"usen \cite{Schmi} has a program to compute the Veech group of
any given square-tiled surface.
She also found positive genus discs as well as noncongruence Veech
groups.
M\"oller \cite{Moe} computes algebraic equations of some square-tiled
surfaces and of their Teich\-m\"uller curves.

\subsection{Acknowledgements}

We thank Anton Zorich for stating questions and some conjectures.
We thank the Institut de Math\'ematiques de Luminy and the
Max-Planck-Institut f\"ur Mathematik for excellent welcome and working
conditions.
We thank Jo\"el Rivat, Martin Schmoll and other participants of the
conference `Dynamique dans l'espace de Teich\-m\"uller et applications
aux billards rationnels' at CIRM in 2003.
We thank Martin M\"oller and Gabriela Schmith\"usen for comments
on a previous version of this paper, circulated under the title
``Square-tiled surfaces in $\cH(2)$''.

\section{Background}

\subsection{Translation surfaces, Veech surfaces}

Let $S$ be an oriented compact surface of genus $g$.
A translation structure on $S$ consists in a set of points
$\{P_1,\ldots,P_n\}$ and a maximal atlas on
$S\raisebox{0.2ex}{$\smallsetminus$} \{P_1,\ldots,P_n\}$ with
translation transition functions.

A holomorphic $1$-form $\omega$ on $S$ induces a translation structure
by considering its natural parameters, and its zeros as points $P_1,
\ldots, P_n$.
All translation structures we consider are induced by holomorphic
$1$-forms.
Slightly abusing vocabulary and notation, we refer to a
translation surface $(S,\omega)$, or sometimes just $S$ or $\omega$.

A translation structure defines: a complex structure, since
translations are conformal; a flat metric with cone-type singularities
of angle $2(k_i+1)\pi$ at order $k_i$ zeros of the $1$-form; and
directional flows $\cF_\theta$ on $S$ for $\theta \in
\left]-\pi,\pi\right]$.

Orbits of the flows $\cF_\theta$ meeting singularities (backward,
resp.\ forward) are called (outgoing, resp.\ incoming) separatrices in
the direction $\theta$.
Orbits meeting singularities both backward and forward are called
\newconcept{saddle connections}; the integrals of $\omega$ along them
are the associated \newconcept{connection vectors}.

Define the singularity type of a $1$-form $\omega$ to be the unordered
tuple $\sigma = (k_1, \ldots, k_n)$ of orders of its zeros (recall
$k_1 + \ldots + k_n = 2 g - 2$, all $k_i > 0$).
The singularity type is invariant by orientation-preserving
diffeomorphisms.
The moduli space $\cH_g$ of holomorphic $1$-forms on $S$ is the
quotient of the set of translation structures by the group
$\Diff^{+}(S)$ of orientation-preserving diffeomorphisms.
$\cH_g$ is stratified by singularity types, the strata are denoted by
$\cH(\sigma)$.

$\SL(2,\R)$ acts on holomorphic $1$-forms: if $\omega$ is a $1$-form,
$\{(U,f)\}$ the translation structure given by its natural parameters,
and $A \in \SL(2,\R)$, then $A\cdot\omega=\{(U, A \circ f)\}$.
As is well known, this action (to the left) commutes with that (to the
right) of $\Diff^{+}(S)$ and preserves singularity types.
Each stratum $\cH(\sigma)$ thus inherits an $\SL(2,\R)$ action.
The dynamical properties of this action have been extensively studied
by Masur and Veech \cite[etc.]{Ma,Ve82}.

From the behavior of the $\SL(2,\R)$-orbit of $\omega$ in
$\cH(\sigma)$ one can deduce properties of directional flows
$\mathcal{F}_\theta$ on the translation surface $(S,\omega)$.
The Veech dichotomy expressed below is a remarkable illustration of
this.

Call \newconcept{affine diffeomorphism} of $(S,\omega)$ an
orientation-preserving homeomorphism $f$ of $S$ such that the
following three conditions hold\\
\textbullet\enspace $f$ keeps the set $\{P_1,\ldots,P_n\}$ invariant;\\
\textbullet\enspace $f$ restricts to a diffeomorphism of
$S\raisebox{0.2ex}{$\smallsetminus$}\{P_1,\ldots,P_n\}$;\\
\textbullet\enspace the derivative of $f$ computed in the natural
charts of $\omega$ is constant.\\
The derivative can then be shown to be an element of $\SL(2,\R)$.

Affine diffeomorphisms of $(S,\omega)$ form its affine group
$\Aff(S,\omega)$, their derivatives form its Veech group
$V(S,\omega)<\SL(2,\R)$, a noncocompact fuchsian group.
The Veech group is the stabilizer of $(S,\omega)$ for the action of
$\SL(2,\R)$ on $\cH_g$.
Veech showed that the derivation map $\Aff(S,\omega) \to V(S,\omega)$
is finite-to-one.
We show (Proposition~\ref{prn:no:translations}) that in $\cH(2)$ it is
one-to-one.

\begin{Veechdichotomy}
If $V(S,\omega)$ is a lattice in $\SL(2,\R)$ (i.e.\
$\vol\left(\lquotient{\SL(2,\R)}{V(S,\omega)}\right)<\infty$) then for
each direction $\theta$, either the flow $\mathcal{F}_\theta$ is
uniquely ergodic, or all orbits of $\cF_\theta$ are compact and $S$
decomposes into a finite number of cylinders of commensurable moduli.
\end{Veechdichotomy}

Cylinder decompositions are further discussed in \S\,\ref{coord:sts}.
Translation surfaces with lattice Veech group are called Veech
surfaces.

\subsection{Square-tiled surfaces, lattice of periods}
\label{sts}
A translation covering is a map $f \colon (S_1,\omega_1)\longrightarrow
(S_2,\omega_2)$ of translation surfaces that\\
\textbullet\enspace is topologically a ramified covering;\\
\textbullet\enspace maps zeros of $\omega_1$ to zeros of $\omega_2$;\\
\textbullet\enspace is locally a translation in the natural parameters
of $\omega_1$ and $\omega_2$.

Translation covers of the standard torus marked at the origin are the
simplest examples of Veech surfaces.
Such surfaces are tiled by squares.

We call them square-tiled.
The Gutkin--Judge theorem states:

\begin{Theorem*}[Gutkin--Judge]
A translation surface $(S,\omega)$ is square-tiled if and only if its
Veech group $V(S,\omega)$ shares a finite-index subgroup with
$\SL(2,\Z)$.
\end{Theorem*}

Translation surfaces with such (arithmetic) Veech groups have also
been called arithmetic; another name for them is origamis.
A proof of Gutkin and Judge's theorem, very different from the
original, is given in appendix~\ref{app:GuJu}.

The subgroup of $\R^2$ generated by connection vectors is the lattice
of relative periods of $(S,\omega)$, denoted by $\Lambda(\omega)$.

\begin{Lemma}
\label{square:tiled:iff:rank:two:sublattice}
A translation surface $(S,\omega)$ is square-tiled if and only if
$\Lambda(\omega)$ is a rank $2$ sublattice of $\Z^2$.
\end{Lemma}

\begin{proof}
If $(S,\omega)$ is square-tiled, connection vectors are obviously
integer vectors, so they span a sublattice of $\Z^2$.
Conversely, let

\noindent
\begin{tabular}[t]{ccc}
$f\colon(S,\omega)$ & $\rightarrow$ &
$\rquotient{\R^2}{\Lambda(\omega)}$,\\
$z$ & $\mapsto$ & $\int_{z_0}^z \omega \bmod \Lambda(\omega)$,
\end{tabular}
where $z_0$ is a given point of $(S,\omega)$.

The integral is well-defined modulo the lattice of absolute periods;
$f$ is a fortiori well-defined.
Since $f$ is holomorphic and onto, it is a covering.
Since relative periods are integer-valued, it is clear that zeros of
$\omega$ project to the origin.
So, given a point $P\neq0$ on the torus, preimages of $P$ are all
regular points, so $P$ is not a branch point.
Hence the covering is ramified only above the origin.
Composing $f$ with the covering $g \colon
\rquotient{\R^2}{\Lambda(\omega)} \longrightarrow
\rquotient{\R^2}{\Z^2}$, we see $(S,\omega)$ is square-tiled.
\end{proof}

A square-tiled surface $(S,\omega)$ is called primitive if
$\Lambda(\omega) = \Z^2$.

\begin{Lemma}
Let $(S,\omega)$ be an $n$-square-tiled surface of genus $g>1$.
If $n$ is prime then $\Lambda(\omega)=\Z^2$.
\end{Lemma}

\begin{proof}
Lemma~\ref{square:tiled:iff:rank:two:sublattice} shows that
$(S,\omega)$ is a ramified cover of
$\rquotient{\R^2}{\Lambda(\omega)}$.
Let $d$ be the degree of the covering.
Then $n = d \cdot [\Z^2 : \Lambda(\omega)]$.
So obviously if $n$ is prime then $\Lambda(\omega) = \Z^2$.
\end{proof}

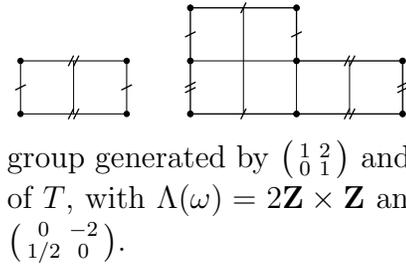
\begin{wrapfigure}{l}{0pt}
  \noindent
  \includegraphics{nonprimitive.eps}
\end{wrapfigure}

Note that $\Lambda(\omega)$ is not always $\Z^2$, as shown by the
examples in the figure.
On the left, a torus $T$ with lattice of periods $2\Z\times\Z$ and
Veech group generated by $\mat{1}{0}{2}{1}$ and $\mat{1}{1/2}{0}{1}$.
On the right, a genus $2$ cover of $T$, with $\Lambda(\omega) =
2\Z\times\Z$ and Veech group generated by $\mat{1}{0}{4}{1}$ and
$\mat{0}{1/2}{-2}{0}$.

\medskip

The following lemma was explained to us first by Martin Schmoll then
by Anton Zorich.

\begin{Lemma}
\label{Schmoll:Zorich}
Let $(S,\omega)$ be a square-tiled surface, then $V(S,\omega)$ is a
subgroup in $V(\rquotient{\R^2}{\Lambda(\omega)},\d z)$.
In particular, if $(S,\omega)$ is primitive, then $V(S,\omega) <
\SL(2,\Z)$.
\end{Lemma}

\begin{proof}
\begin{tabular}[t]{cccc}
Let $\phi\colon$ & $V(S,\omega)$ & $\rightarrow$ &
$V(\rquotient{\R^2}{\Lambda(\omega)},\d z)$,\\
& $A = df,\ f \in \Aff(S,\omega)$ & $\mapsto$ & $A$.
\end{tabular}

The only difficulty is to show that $\phi$ is well-defined i.e.\ that
any element $A$ in $V(S,\omega)$ preserves $\Lambda(\omega)$.
Since any element of the affine group maps a connection to a
connection, hence $A$ maps a connection vector to a connection vector
(i.e.\ an element in $\Lambda(\omega)$).
\end{proof}

\begin{Remark}
As shown by the examples above, there are Veech groups of square-tiled
surfaces which are \emph{not} subgroups of $\SL(2,\Z)$.
\end{Remark}

\subsection{Cylinders of square-tiled surfaces}
\label{coord:sts}
A square-tiled surface decomposes into maximal horizontal cylinders,
bounded above and below by unions of saddle connections, each of which
appears once on the top of a cylinder and once on the bottom of a
cylinder.
Gluing the cylinders alongs these saddle connections builds back the
surface.

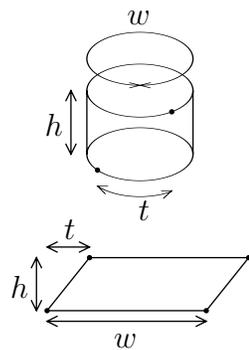
\begin{wrapfigure}{l}{95pt}
  \setlength{\intextsep}{0pt}
  \psfrag{w}{$w$}
  \psfrag{h}{$h$}
  \psfrag{t}{$t$}
  \includegraphics{cyldim.eps}\\
  \includegraphics{cyldimflat.eps}
\end{wrapfigure}

A cylinder on a translation surface is isometric to
$\rquotient{\R}{w\Z} \times [0,h]$, for some $h$ and $w$.

\textbf{Convention}.
We refer to these dimensions as height and width respectively,
whether the `horizontal direction of the cylinder' coincides with the
horizontal direction of the surface or not.

An additional twist parameter $t$ is needed, measuring the distance
along the `horizontal direction of the cylinder' between some
(arbitrary) reference points on the bottom and top of the cylinder,
for instance some ends of saddle connections.

\subsection{Action of $\SL(2,\Z)$ on square-tiled surfaces}

\begin{Lemma}
\label{primitive:orbit}
The $\SL(2,\Z)$-orbit of a primitive $n$-square-tiled surface is the
set of primitive $n$-square-tiled surfaces in its $\SL(2,\R)$-orbit.
\end{Lemma}

\begin{proof}
$\SL(2,\Z)$ preserves $\Z^2$ ($= \Lambda(\omega)$ if $(S,\omega)$ is
primitive square-tiled) and hence the property of being primitive
square-tiled.
Conversely, if $(S,\omega)$ is primitive square-tiled and
$(S_1,\omega_1)=A{\cdot}(S,\omega)$ is square-tiled for some $A \in
\SL(2,\R)$, then $\Lambda(\omega_1) = A{\cdot} \Lambda(\omega)$ means
$A$ preserves $\Z^2$, so $A \in \SL(2,\Z)$.
\end{proof}

\begin{Remark}
The number of squares, $n$, is preserved by $\SL(2,\R)$ because it is
the area of the surface.
Consequently $\SL(2,\Z)\cdot(S,\omega)$ is finite.
\end{Remark}

\textbf{Notation}.
Denote by $\U=\mat{1}{0}{1}{1}$ and $\V=\mat{0}{1}{-1}{0}$ the
standard generators of $\SL(2,\Z)$, and by $\cU = \langle \U \rangle =
\{\mat{1}{0}{n}{1}\colon n \in \Z\}$ the subgroup generated by $\U$.

\setlength{\intextsep}{-5pt}

\begin{wrapfigure}{r}{0pt}
  \psfrag{u}{{\large $\xrightarrow{\U}$}}
  \includegraphics{actionsquareu.eps}
  \hfil
  \psfrag{v}{{\large $\xrightarrow{\V}$}}
  \includegraphics{actionsquarev.eps}
\end{wrapfigure}

Here is the action of $\U$ and $\V$ on squares.

\setlength{\intextsep}{0pt}

The action on square-tiled surfaces is obtained by applying the same
to all square tiles.
The new horizontal cylinder decomposition is then recovered by cutting
and gluing (see example in \S\,\ref{sec:tools:action}).

\subsection{Hyperelliptic surfaces, Weierstrass points}
\label{hyperell}

Recall that a Riemann surface $X$ of genus $g$ is hyperelliptic
if there exists a degree $2$ meromorphic function on $X$.
Such a function induces a holomorphic involution on $X$.
This involution has $2g+2$ fixed points called Weierstrass points.
The set of these points is invariant by all automorphisms of the
complex structure.
A translation surface is called hyperelliptic if the underlying
Riemann surface is hyperelliptic.

Hyperelliptic translation surfaces have been studied by Veech.
He showed \cite{Ve95} that in genus $g$ they are obtained from
centrosymmetric polygons with $4g$ or $4g+2$ sides by pairwise
identifying opposite sides.

The hyperelliptic involution is in these coordinates the reflection in
the center of the polygon; the Weierstrass points are the center of
the polygon, the midpoints of its sides, and the vertices (identified
into one point) in the $4g$ case (in the $4g+2$ case the vertices are
indentified into two points exchanged by the hyperelliptic
involution).

\subsection{Cusps}

Let $\Gamma$ be a fuchsian group.
A parabolic element of $\Gamma$ is a matrix of trace $2$ (or $-2$).
A point of the boundary at infinity of $\HH^2$ is parabolic if it is
fixed by a parabolic element of $\Gamma$.
A cusp is a conjugacy class under $\Gamma$ of primitive parabolic
elements (primitive meaning not powers of other parabolic elements of
$\Gamma$).

Recall that a lattice admits only a finite number of cusps.

Geometrically, each cusp in $\lquotient{\HH^2}{\Gamma}$ has, for some
positive $\lambda$ called its \newconcept{width}, neighborhoods
isometric to the quotients of the strips $\{z\in \C\colon 0<|\Re
z|<\lambda,\ \Im z > M\}$ by the translation $z \mapsto z + \lambda$,
for large $M$.

On a Veech surface $(S,\omega)$, any `periodic' direction is fixed by
a parabolic element of the Veech group.
Conversely the eigendirection of a parabolic element in the Veech
group is a `periodic' direction.
We call such directions parabolic.
Thus parabolic limit points of $V(S,\omega)$ are cotangents of
parabolic directions.

When $(S,\omega)$ is a square-tiled surface, the set of parabolic
limit points is $\Q$.
Cusps are therefore equivalence classes of rationals under the
homographic action of $V(S,\omega)$.
The following lemma gives a combinatorial description of cusps for a
square-tiled surface.

\begin{Lemma}[Zorich]
\label{cusp:lemma}
Let $(S,\omega)$ be a primitive $n$-square-tiled surface and
$E=\SL(2,\Z)\cdot (S,\omega)$ the set of $n$-square-tiled surfaces in
its orbit.
The cusps of $(S,\omega)$ are in bijection with the $\cU$-orbits of
$E$.
\end{Lemma}

\begin{proof}
Denote by $\cC$ the set of cusps of $(S,\omega)$.

Let $\varphi\colon$
\begin{tabular}[t]{ccccc}
  $\SL(2,\Z)$
  & $\xrightarrow{f}$
  & $\Q$
  & $\xrightarrow{\pi}$
  & $\cC$,
  \\
  $A$
  & $\mapsto$
  & $A^{-1} \infty$
  & $\mapsto$
  & $A^{-1} \infty \bmod V(S,\omega)$.
\end{tabular}

Note that $\infty$ corresponds to the horizontal direction in
$(S,\omega)$ because the projective action is the action on co-slopes
and not on slopes.
$A^{-1} \infty$ corresponds to the direction on $(S,\omega)$ that is
mapped by $A$ to the horizontal direction of $A \cdot (S,\omega)$.

$\varphi$ pulls down as $\psi\colon$
\begin{tabular}[t]{ccc}
  $E$ & $\rightarrow$ & $\cC$,\\
  $A\cdot (S,\omega)$ & $\mapsto$ & $A^{-1} \infty \bmod V(S,\omega)$.
\end{tabular}

$\psi$ is well-defined: if $A \cdot (S,\omega) = B \cdot (S,\omega)$,
then $\exists P \in V(S)$, $B = A P$, so setting $A^{-1} \infty =
\alpha$, $B^{-1} \infty = \beta$, we have $\beta = B^{-1} \infty =
(B^{-1}A)A^{-1} \infty = P^{-1} \alpha$, so $\alpha$ and $\beta$
correspond to the same cusp.
Further, $\psi$ is surjective because $f$ is.
Indeed, $\forall \alpha = p/q$, $\exists A \in \SL(2,\Z)$ s.t.\
$A^{-1} \infty = \alpha$.
(The orbit of $\infty$ under $\SL(2,\Z)$ is $\Q$.)

Recall that the stabilizer of $\infty$ for the action of $\SL(2,\Z)$
is $\cU$.
If $\psi(S_1,\omega_1) = \psi(S_2,\omega_2)$, where $(S_1,\omega_1) =
A \cdot (S,\omega)$ and $(S_2,\omega_2) = B \cdot (S,\omega)$, then
$\varphi(A)=\varphi(B)$.

Let $\alpha = f(A) = A^{-1} \infty$ and $\beta = f(B) = B^{-1}
\infty$.
Since $\alpha$ and $\beta$ correspond to the same cusp, $\exists P \in
V(S)$ s.t.\ $\beta = P\alpha$.
So $\infty = A \alpha = A P^{-1} \beta = A P^{-1} B^{-1} \infty$ which
implies $A P^{-1} B^{-1} \in \cU$ i.e.\ $\exists U^k \in \cU$ s.t.\ $A
P^{-1} = U^k B$ i.e.\ $A P^{-1} \cdot (S,\omega) = A \cdot (S,\omega)
= U^kB \cdot (S,\omega)$, so that $(S_1,\omega_1)$ and
$(S_2,\omega_2)$ are in the same $\cU$-orbit.

Conversely: if $(S_2,\omega_2) =U^k(S_1,\omega_1)$ with $U^k \in
\cU$, and $(S_2,\omega_2)=B\cdot (S,\omega)$ and $(S_1,\omega_1)=A
\cdot (S,\omega)$, then $\psi(S_2,\omega_2) = B^{-1} \infty =
A^{-1}U^{-k}\infty = A^{-1}\infty = \psi (S_1,\omega_1)$.
\end{proof}

\subsection{Elliptic points}

Recall that in a fuchsian group $\Gamma$, any elliptic element has
finite order and is conjugate to a rational rotation.

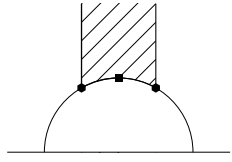
\begin{wrapfigure}{l}{0pt}
  \noindent
  \includegraphics{modular.eps}
\end{wrapfigure}

A fixed point in $\HH^2$ of an elliptic element of $\Gamma$ is called
elliptic.
Its projection to the quotient $\lquotient{\HH^2}{\Gamma}$ is a cone
point, with a curvature default.
For instance the modular surface $\lquotient{\HH^2}{\SL(2,\Z)}$ has
two cone points, of angles $\pi$ and $2\pi/3$.

Suppose that $\Gamma$ is the Veech group of a translation surface and
has an elliptic point.
By applying a convenient element of $\SL(2,\R)$, we can suppose that
this point is $i$.
The corresponding elliptic element is a rational rotation.
The translation surfaces which project to $i$ have this rotation in
their Veech group.
This roughly means that they have an apparent symmetry.
At the Riemann surface level, the rotation is an automorphism of the
complex structure (it modifies the vertical direction but not the
metric).
For genus $1$, the cone point $i$ (resp.\ $e^{i\pi/3}$) of the modular
surface corresponds to the square (resp.\ hexagonal) torus, which has
a symmetry of projective order $2$ (resp.\ $3$).

One should note that the translation surfaces obtained from rational
polygonal billiards always have elliptic elements in their Veech
group:
writing the angles of a simple polygon as $(k_1\pi/r, \ldots,
k_q\pi/r)$, with $k_1$, \ldots, $k_q$, $r$ coprime,
the covering translation surface is obtained by gluing $2r$ copies by
symmetry.
The rotation of angle $2\pi/r$ is in the Veech group (this rotation is
minus the identity if $r=2$).
Many explicit calculations of lattice Veech groups make use of this
remark (see \cite{Ve89}, \cite{Vo}, \cite{Wa}).
Our method is completely different.

\subsection{The Gauss--Bonnet Formula}

Let $\Gamma$ be a finite-index subgroup of $\SL(2,\Z)$ containing
$-\Id$.
The quotient of $\lquotient{\SL(2,\R)}{\Gamma}$ is the unit tangent
bundle to an orbifold surface with cusps $S_\Gamma$.
Algebraic information on the group is related to the geometry of the
surface.

Let $d$ be the index $[\SL(2,\Z):\Gamma]$ of $\Gamma$ in $\SL(2,\Z)$,
$e_2$ (resp.\ $e_3$) the number of conjugacy classes of elliptic
elements of order $2$ (resp.\ $3$) of $\Gamma$, $e_\infty$ the number
of conjugacy classes of cusps of $\Gamma$.

Then the surface $S_\Gamma$ has hyperbolic area $d\frac{\pi}{3}$,
$e_2$ cone points of angle $\pi$, $e_3$ cone points of angle
$\frac{2\pi}{3}$, $e_\infty$ cusps, and its genus $g$ is given by:

\begin{GBFormula}
  $g = 1 + d/12 - e_2/4 - e_3/3 - e_\infty/2$.
\end{GBFormula}

\section{Specific Tools}

In this section we give specific properties of the stratum $\cH(2)$,
and a combinatorial coordinate system for square-tiled surfaces in
$\cH(2)$.

\subsection{Hyperellipticity}

First recall that any genus $2$ Riemann surface is hyperelliptic.
Given a genus $2$ Riemann surface $X$ and its hyperelliptic involution
$\tau$, any $1$-form $\omega$ on $X$ satisfies $\tau^*\omega = -
\omega$.

In the moduli space of holomorphic $1$-forms of genus $2$, $\cH(2)$ is
the stratum of $1$-forms with a degree $2$ zero (a cone point of angle
$6\pi$).

As said in \S\,\ref{hyperell}, any translation surface in $\cH(2)$
can be represented as a centro-symmetric octagon.
The six Weierstrass points are the center of the polygon, the middles
of the sides and the cone-type singularity.
The position of the Weierstrass points in a surface decomposed into
horizontal cylinders is described in \S\,\ref{sec:twovalinv}.

\subsection{Separatrix diagrams}
\label{section:separatrix:diagrams}

Forms in $\cH(2)$ have a single degree $2$ zero, geometrically a cone
point of angle $6\pi$, with three outgoing separatrices and three
incoming ones in any direction.

Recall that the horizontal direction of a square-tiled surface is
completely periodic; the horizontal separatrices are saddle
connections.
The combinatorics of these connections is called a separatrix diagram
in \cite{KoZo}.
The surface is obtained from this diagram by gluing cylinders along
the saddle connections.

Each outgoing horizontal separatrix returns to the saddle making an
angle $\pi$, $3\pi$ or $5\pi$ with itself.
Four separatrix diagrams are combinatorially possible (up to rotation
by $2\pi$ around the cone point); they correspond to return angles
$(\pi, \pi, \pi)$, $(\pi, 3\pi, 5\pi)$, $(3\pi, 3\pi, 3\pi)$, $(5\pi,
5\pi, 5\pi)$:

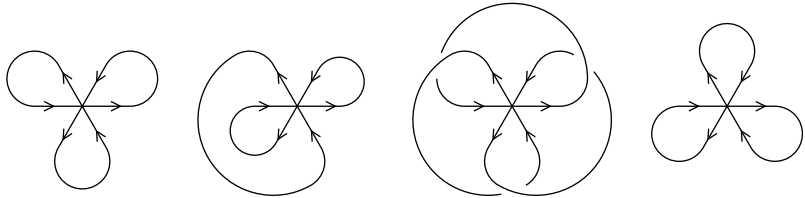
\begin{center}
\includegraphics{separatrixdiagrams3.eps}
\end{center}

There is no consistent way of gluing cylinders along the saddle
connections of the first and last diagrams to obtain a translation
surface.

The second diagram is possible with the condition that the saddle
connections that return with angles $\pi$ and $5\pi$ have the same
length; this diagram corresponds to surfaces with two cylinders.
The third diagram corresponds to surfaces with one cylinder, with no
restriction on the lengths of the saddle connections.

\subsection{Parameters for square-tiled surfaces in $\cH(2)$}

Here we give complete combinatorial coordinates for square-tiled
surfaces in $\cH(2)$.
See figures in \S\,\ref{sec:twovalinv}.

\medskip

\textbf{Notation.}
We use $\wedge$ for greatest common
divisor, and $\vee$ for least common multiple.

\subsubsection{One-cylinder surfaces}

A one-cylinder surface is parametrized by the height of the cylinder,
the lengths of the three horizontal saddle connections (a triple of
integers up to cyclic permutation), and the twist parameter.
If all three horizontal saddle connections have the same length, the
twist parameter is taken to be less than that length; otherwise, less
than the sum of the three lengths.

For primitive surfaces, the height is $1$, and the lengths of the
three horizontal saddle connections add up to the area $n$ of the
surface.

The horizontal saddle connections appear in some (cyclic) order on the
bottom of the cylinder, and in reverse order on the top.

\subsubsection{Two-cylinder surfaces}
\label{param:two:cyl:surf}

Labeling the horizontal saddle connections according to their return
angles, call them $\gamma_\pi$, $\gamma_{3\pi}$, $\gamma_{5\pi}$.
Call $\ell_1$ the common length of $\gamma_{\pi}$ and $\gamma_{5\pi}$,
and $\ell_2$ the length of $\gamma_{3\pi}$.
One cylinder is bounded below by $\gamma_{\pi}$ and above by
$\gamma_{5\pi}$; the other one is bounded below by $\gamma_{5\pi}$ and
$\gamma_{3\pi}$, and above by $\gamma_{\pi}$ and $\gamma_{3\pi}$.

A two-cylinder surface is determined by the heights $h_1$, $h_2$ and
widths $w_1=\ell_1$, $w_2=\ell_1+\ell_2>w_1$ of the cylinders as well
as two twist parameters $t_1$, $t_2$ satisfying $0 \leqslant t_1 <
w_1$, $0 \leqslant t_2 < w_2$.
The area of the surface is $h_1 w_1 + h_2 w_2 = n$.
For primitive surfaces, $h_1\wedge h_2 = 1$.
For prime $n$, in addition, $\ell_1 \wedge \ell_2 = 1$, and (P)
$\ell_1Ê\wedge h_2 = 1$.

\subsection{Action of $\SL(2,\Z)$}
\label{sec:tools:action}

The action of $\V$ (rotation by $\pi/2$) does not preserve separatrix
diagrams in general.
The horizontal cylinder decomposition of $\V\cdot S$ is the vertical
cylinder decomposition of $S$.

$\U$ is the primitive parabolic element in $\SL(2,\Z)$ that preserves
the horizontal direction.
Its action preserves separatrix diagrams, as well as heights $h_i$ and
widths $w_i$ of horizontal cylinders $C_i$, and only changes twist
parameters $t_i$ to $(t_i + h_i) \bmod w_i$.

Here is an example of how $\U$ acts on a surface.

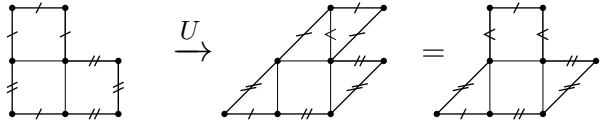
\begin{center}
  \psfrag{a}{{\large $\xrightarrow{\U}$}}
  \psfrag{b}{$=$}
  \includegraphics{action2cyl.eps}
\end{center}

For prime $n$, given a cyclically ordered $3$-partition $(a,b,c)$ of
$n$, all one-cylinder surfaces with bottom sides of lengths $a$, $b$,
$c$ (up to cyclic permutation) are in the same $\cU$-orbit, or cusp
(see Lemma~\ref{cusp:lemma}).

The following lemma describes $\U$-orbits of two-cylinder surfaces in
$\cH(2)$ by giving their sizes and canonical representatives.

\begin{Lemma}
\label{lem:cusp:widths}
Let $S$ be a primitive two-cylinder $n$-square-tiled surface in
$\cH(2)$ with parameters $h_i$, $w_i$, $t_i$ ($i=1,2$).
Then the cardinality of its $\U$-orbit (its \newconcept{cusp width})
is

\smallskip

\noindent
\centerline{%
$\displaystyle \mathrm{cw}(S) =
\frac{w_1}{w_1 \wedge h_1} \vee \frac{w_2}{w_2 \wedge h_2}$
\quad ($=\frac{w_1}{w_1 \wedge h_1} \times \frac{w_2}{w_2 \wedge h_2}$
for prime $n$).
}

\smallskip

The surface $S'$ with $h'_i = h_i$, $w'_i = w_i$, and $t'_i = t_i
\bmod (w_i \wedge h_i)$ is a ``canonical'' representative of the
$\U$-orbit of $S$.
Each surface thus has a unique representative with $0 \leqslant t'_i <
w_i \wedge h_i$.
\end{Lemma}

\begin{proof}
Observe that $\U^k\cdot S$ has widths $w_i$, heights $h_i$, and twist
parameters $(t_i + k h_i) \bmod w_i$.
So for $U^k\cdot S$ to coincide with $S$, the integer $k$ must be a
multiple of $\frac{w_i}{w_i \wedge h_i}$ for each $i$.
The cusp width is the least such positive $k$, the least common
multiple of $\frac{w_1}{w_1 \wedge h_1}$ and $\frac{w_2}{w_2 \wedge
h_2}$.
The second part is a simple application of the Chinese remainder
theorem.
\end{proof}

\section{Results}
\label{sec:results}
This section expands the results summarized in the introduction,
detailed proofs are postponed to the next sections.
Additional conjectures appear in \S\,\ref{section:numerical:evidence}.

\subsection{Two orbits}

Theorem~\ref{thm:twoorbits} can be reformulated as:

\begin{Proposition}
\label{prn:twoorbits}
Given a prime $n \geqslant 5$, the primitive $n$-square-tiled surfaces
in $\cH(2)$ fall into two $\SL(2,\Z)$ orbits.
\end{Proposition}

The idea for proving this is first to give an invariant which takes
two different values, thus proving that there are at least two orbits
(see \S\,\ref{sec:invariant} below, and \S\,\ref{sec:twovalinv}), then
prove that there are exactly two orbits by showing that each orbit
contains a one-cylinder surface (see \S\,\ref{sec:redonecyl}), and
that all one-cylinder surfaces with the same invariant are indeed in
the same orbit (\S\,\ref{sec:linkonecyl}).

We will call these orbits A and B.

\begin{Remark}
An extension of this result in some components of higher-dimensional
strata is presented in appendix~\ref{app:hypel}.
\end{Remark}

\subsection{Invariant}
\label{sec:invariant}
We present a geometric invariant that can easily be computed for any
primitive square-tiled surface in $\cH(2)$ (for instance presented in
its decomposition into horizontal cylinders.)

The Weierstrass points of a surface in $\cH(2)$ are\\
\textbullet\enspace the saddle ($6 \pi$-angle cone point),\\
\textbullet\enspace  and five regular points.

\begin{Lemma}
\label{lemma:invariant}
The number of integer Weierstrass points of a primitive square-tiled
surface is invariant under the action of $\SL(2,\Z)$.
\end{Lemma}

By integer point we mean a vertex of the square tiling.
The proof of the lemma is obvious, since $\SL(2,\Z)$ preserves $\Z^2$.

\begin{Proposition}
\label{prn:invariant}
Primitive $n$-square-tiled surfaces in $\cH(2)$ have\\
\textbullet\enspace for $n=3$, exactly $1$ integer Weierstrass
point,\\
\textbullet\enspace for even $n$, exactly $2$,\\
\textbullet\enspace for odd $n$, either $1$ or $3$ (both values
occur).
\end{Proposition}

Martin M\"oller pointed out to us that this invariant also appears in
\cite[\S\,2, formula (6)]{Ka} in algebraic geometric language; Kani's
normalized covers correspond to our orbit B.
This invariant is also mentioned in \cite[Remark~3.4]{Moe}.

\subsection{Elliptic affine diffeomorphisms}
\label{elliptic}

\begin{Proposition}
\label{prn:no:translations}
A translation surface in $\cH(2)$ has no nontrivial translation in its
affine group.
Hence the derivation from its affine group to its Veech group is an
isomorphism.
\end{Proposition}

\begin{Proposition}
\label{prn:ellipticthree}
A translation surface in $\cH(2)$ can have no elliptic element of
order $3$ in its Veech group.
\end{Proposition}

\begin{Lemma}
\label{lemma:octagon}
Any $\V$-invariant Veech surface in $\cH(2)$ can be represented as a
$\V$-invariant octagon.
\end{Lemma}

\begin{Proposition}
\label{prn:elliptictwo}
For any given prime $n$, there exist $\V$-invariant $n$-square-tiled
$\cH(2)$ surfaces.
All of them have the same invariant, namely, A if $n \equiv - 1 \pmod
4$ and B if $n \equiv 1 \pmod 4$.
\end{Proposition}

\begin{Remark}
This proposition implies the following interesting fact: there are
finite-covolume Teich\-m\"uller discs with no elliptic points.
This differs from the billiard case which has been the main source of
explicit examples of lattice Veech groups.
\end{Remark}

\subsection{Countings}

The asymptotic number of square-tiled surfaces in $\cH(2)$ of area
bounded by $N$ is given in \cite{Zo} (see also \cite{EsOk} and
\cite{EsMaSc}) to be $\zeta(4) \frac{N^4}{24}$ for one-cylinder
surfaces and $\frac{5}{4} \zeta(4) \frac{N^4}{24}$ for two-cylinder
surfaces.
The mean order for the number of square-tiled surfaces of area exactly
$n$ is therefore $\zeta(4) \frac{n^3}{6}$ for one-cylinder surfaces
and $\frac{5}{4} \zeta(4) \frac{n^3}{6}$ for two-cylinder surfaces.

The following proposition, from which Theorem~\ref{thm:countings}
follows, states that for prime $n$, there are in fact asymptotics for
these numbers, which are $\zeta(4)$ times smaller than the mean order.

\begin{Proposition}
\label{prn:countings}
For prime $n$, there are $O(n)$ elliptic points, and the following
countings and asymptotics hold for surfaces and cusps, according to
the number of cylinders and to the orbit.
\begin{center}
\noindent
\begin{tabular}{c@{\,}c}
  surfaces: & cusps: \\
  \noindent
  \begin{tabular}{|c|c|c|c|}
    \hline
    & $1$-cyl & $2$-cyl & all \\
    \hline
    A & $\frac{n(n-1)(n+1)}{24}$
    & $\sim\frac{7}{8}\frac{n^3}{6}$
    & $\sim\frac{9}{8}\frac{n^3}{6}$ \\[1pt]
    \hline
    B & $\frac{n(n-1)(n-3)}{8}$ &
    $\sim\frac{3}{8}\frac{n^3}{6}$
    & $\sim\frac{9}{8}\frac{n^3}{6}$ \\[1pt]
    \hline
    all & $\frac{n(n-1)(n-2)}{6}$ &
    $\sim\frac{5}{4}\frac{n^3}{6}$
    & $\sim\frac{9}{4}\frac{n^3}{6}$ \\[1pt]
    \hline
  \end{tabular}
  &
  \noindent
  \begin{tabular}{|c|c|c|c|}
    \hline
    & $1$-cyl & $2$-cyl & all \\
    \hline
    A & $\frac{(n-1)(n+1)}{24}$
    & $o(n^{3/2 + \varepsilon})$
    & $\sim \frac{n^2}{24}$ \\[1pt]
    \hline
    B & $\frac{(n-1)(n-3)}{8}$
    & $o(n^{3/2 + \varepsilon})$
    & $\sim \frac{n^2}{8}$ \\[1pt]
    \hline
    all & $\frac{(n-1)(n-2)}{6}$
    & $o(n^{3/2 + \varepsilon})$
    & $\sim \frac{n^2}{6}$ \\[1pt]
    \hline
  \end{tabular}
\end{tabular}
\end{center}
\end{Proposition}

This proposition gives more detail than Theorem~\ref{thm:countings} by
distinguishing one-cylinder and two-cylinder cusps and surfaces.
Proposition~\ref{prn:one:cyl:dir} and
Proposition~\ref{prn:dev:mean:order} are corollaries of this
proposition.

\begin{Remarks}
Orbits A and B have asymptotically the same size (same number of
square-tiled surfaces).
However orbit B has asymptotically three times as many one-cylinder
surfaces as orbit A.

In each orbit the proportion of two-cylinder cusps is asymptotically
negligible; however it is not the case for the proportion of
two-cylinder surfaces.
This shows that the average width of the two-cylinder cusps grows
faster than $n$.
(One-cylinder cusps all have width $n$.)
\end{Remarks}

\section{Proof of main theorem (two orbits)}

In this section we first prove Proposition~\ref{prn:invariant}, then
Proposition~\ref{prn:twoorbits}.

\textbf{Convention for figures}.
In all figures, we represent a square-tiled surface $S$ in $\cH(2)$ by
a fundamental octagonal domain.
$S$ is obtained by identifying pairs of parallel sides of same length;
all vertices (black dots) get identified to the saddle.
Circles are sometimes used to indicate the other Weierstrass points.

\begin{wrapfigure}{r}{150pt}
  \includegraphics{wpts1.eps}\\
  \includegraphics{wpts2.eps}
\end{wrapfigure}

Except in \S\,\ref{square:octagon}, the octagon reflects
horizontal cylinders: nonhorizontal sides are identified by
horizontal translations.
On one-cylinder surfaces, the horizontal sides on the top and on the
bottom of the cylinder are identified in opposite cyclic order.
Two-cylinder surfaces are represented with the cylinder of least width
on top of the widest one, to the left.
Its top side is glued to the leftmost side under the bottom cylinder.
The remaining two sides, to the right on the top and bottom of the
bottom cylinder, are identified with each other.

\subsection{Two values of the invariant}
\label{sec:twovalinv}

Here we prove Proposition~\ref{prn:invariant}, about the possible
values of the number of integer Weierstrass points of a primitive
square-tiled surface in $\cH(2)$.

Recall that the hyperelliptic involution turns the cylinders
upside-down.
We deduce the position of Weierstrass points (see figure).

The saddle is always an integer Weierstrass point.
We discuss the case of the remaining five, depending on the parity of
the parameters.

Under the hyperelliptic involution:\\
\textbullet\enspace saddle connections that bound a cylinder both on
its top and on its bottom are mapped to themselves with reversed
orientation, so that their middpoint is fixed: it is a Weierstrass
point, integer when the length of the saddle connection is even.\\
\textbullet\enspace the core circle of a cylinder, also mapped to itself
with orientation reversed, has two antipodal fixed points.
If the cylinder has odd height, none of them is integer.
When the height is even and the width odd, one of them is integer.
When the height and width are even, either both or none is integer,
depending on the parity of the twist parameter.

\subsubsection{One-cylinder case}

The core of the (height $1$) cylinder contains two non-integer
Weierstrass points.
The remaining three are the midpoints of the horizontal connections
(whose lengths add up to $n$).

If $n$ is odd, it splits into either $3$ odd lengths (no integer
Weierstrass point), or $1$ odd and $2$ even lengths ($2$ integer
Weierstrass points).
For $n=3$ all lengths are $1$ (hence odd); for greater odd $n$ both
cases occur.

If $n$ is even, two lengths are odd and one even (if all were even,
the surface could not be primitive).
This completes the one-cylinder case.

\subsubsection{Two-cylinder case}

We use parameters $h_1$, $h_2$, $w_1$, $w_2$, $t_1$, $t_2$ introduced
above.
We also use $\ell_1$ and $\ell_2$ to denote the lengths of the
horizontal saddle connections.
We then have:

\centerline{%
$\ell_1 = w_1$, $\ell_1 + \ell_2 = w_2$, $n = w_1 h_1 + w_2 h_2 = h_1
\ell_1 + h_2 (\ell_1 + \ell_2)$ ($*$).}

\textbullet\enspace \textbf{Odd $n$.}
If $\ell_2$ is even, the corresponding Weierstrass point is integer.
Because $n$ is odd, equation ($*$) implies that $\ell_1$ is odd, thus
both cylinders have odd widths, and still by ($*$) one of the heights
must be even.
The corresponding cylinder has one integer Weierstrass point on its
core line.
The total number of integer Weierstrass points is then $3$.

If $\ell_2$ is odd, the corresponding Weierstrass point is
non-integer; if $\ell_1$ is odd (resp.\ even), then $w_2$ is even
(resp.\ odd), thus by ($*$) $h_1$ (resp.\ $h_2$) has to be odd,
meaning the top (resp.\ bottom) cylinder contains two non-integer
Weierstrass points.
The two Weierstrass points in the bottom (resp.\ top) cylinder are
integer if $h_2$ is even and $t_2$ is odd (resp.\ if $h_1$ and $t_1$
are even), non-integer otherwise (see figure above).
The value of the invariant is accordingly $3$ or $1$.

For $n=3$, $\ell_1=\ell_2=1$; for greater odd $n$ both values do
occur.

\textbullet\enspace \textbf{Even $n$.}
Recall that primitivity implies $h_1\wedge h_2 = 1$.
In particular at least one of them is odd.

If both heights are odd, the Weierstrass points inside the cylinders
are non-integer, and because $n = (h_1 + h_2) \ell_1 + h_2 \ell_2$ is
even, $\ell_2$ has to be even, so the last Weierstrass point is
integer, and the invariant is $2$.

If $h_1$ is odd and $h_2$ even, then, by ($*$), $\ell_1$ has to be
even.
Then if $\ell_2$ is odd, the corresponding Weierstrass point is
non-integer, one of the Weierstrass points inside the bottom cylinder
is integer, and the invariant is $2$.
If $\ell_2$ is even, the corresponding Weierstrass point is integer,
and $t_2$ has to be odd for the surface to be primitive, hence the
remaining Weierstrass points are non-integer, and the invariant is
$2$.

The last case to consider is when $h_1$ is even and $h_2$ odd.
If $\ell_1$ is odd, then so is $\ell_2$ (by ($*$)), so one Weierstrass
point in the top cylinder is integer, and the invariant is $2$.
If $\ell_1$ is even, then $\ell_2$ is also even by ($*$).
The corresponding Weierstrass point is integer, and $t_1$ is odd for
primitiveness.
Thus all Weierstrass points inside cylinders are non-integer, and the
invariant is $2$.

This completes the two-cylinder case, and
Proposition~\ref{prn:invariant} is proved.

\textbullet\enspace \textbf{Summary of two-cylinder case.}
For future reference, we sum up the case study above in a table giving
the invariant for odd $n$ according to the parity of $h_1$, $h_2$,
$\ell_1$, $\ell_2$ (recall that $w_1=\ell_1$ and $w_2=\ell_1+\ell_2$).

\setlength{\intextsep}{2pt}
\begin{wraptable}{l}{0pt}
\noindent
\begin{tabular}{ccccc}
  \hline
  $h_1$ & $h_2$ & $\ell_1$ & $\ell_2$ & invariant\\
  \hline
  0 & 1 & 1 & 0 & 3 \\
  1 & 0 & 1 & 0 & 3 \\
  0 & 1 & 0 & 1 & $t_1$ odd: 1; $t_1$ even: 3 \\
  1 & 0 & 1 & 1 & $t_2$ odd: 3; $t_2$ even: 1 \\
  1 & 1 & 0 & 1 & 1 \\
  1 & 1 & 1 & 1 & 1 \\
  \hline
\end{tabular}
\end{wraptable}

Table for odd $n$ case.

The other combinations of parities of the parameters cannot happen
for odd $n$ and primitive surfaces.

Note that for even $n$ we concluded that the invariant is $2$ for all
primitive surfaces.

\setlength{\intextsep}{0pt}

\subsection{Reduction to one cylinder}
\label{sec:redonecyl}

\begin{Proposition}
Each orbit contains a one-cylinder surface.\\
Equivalently, each surface has a direction in which it decomposes in
one single cylinder.
\end{Proposition}

A baby version of this proposition is the following lemma.
\begin{Lemma}
A two-cylinder surface of height $2$ tiled by a prime number of
squares has one-cylinder directions.
\end{Lemma}

\begin{proof}[Proof of the lemma]
Consider a surface made of two cylinders, both of height $1$.
Since $n$ is prime, the two widths are relatively prime.
By acting by $\U$, the twists can be set to any values (see
Lemma~\ref{lem:cusp:widths}).
Set the top twist to $0$ and the bottom twist to $1$.
Then by considering the vertical flow, we get a one-cylinder surface.
\end{proof}

We prove the proposition by induction on the height of the surface:
given a two-cylinder surface, we show that its orbit contains a
surface of strictly smaller height.

Consider a two-cylinder square-tiled surface $S$ in $\cH(2)$, with a
prime number of square tiles.
By acting by $\U$ we can move to the canonical representative of the
same cusp (see Lemma~\ref{lem:cusp:widths}), so we will assume
$t_i<w_i$, $i=1,2$.

We split our study into four cases according to which twists are zero.

Case 1.
Both twists are nonzero.

\begin{center}
  \psfrag{A}{{\Small $A$}}
  \psfrag{B}{{\Small $B$}}
  \psfrag{C}{{\Small $C$}}
  \psfrag{D}{{\Small $D$}}
  \psfrag{E}{{\Small $E$}}
  \psfrag{F}{{\Small $F$}}
  \psfrag{G}{{\Small $G$}}
  \psfrag{H}{{\Small $H$}}
  \includegraphics{twoorbits1.eps}
\end{center}

Call $h_1$, $h_2$ the heights and $t_1$, $t_2$ the twists of the
horizontal cylinders of $S$.
Consider the rotated surface $\V S$.
If $\V S$ consists of one horizontal cylinder, we are done.
Otherwise, it has two horizontal cylinders, which are the vertical
cylinders of $S$, and fill $S$.
Looking to the right of $A$, $H$, and $B$, we see all vertical
cylinders of $S$.
The vertical cylinder to the right of $A$ has height at most $t_2$, that
to the right of $B$ also, and that to the right of $H$ at most $t_1$.
So one of the vertical cylinders has heights at most $t_2$, and the
other one has height at most $t_1$.
The sum of their heights is hence at most $t_1 + t_2$, so it is less
than $h_1 + h_2$.

Case 2.
The bottom twist is nonzero but the top twist is zero.

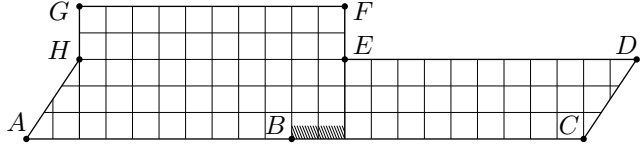
\begin{center}
  \psfrag{A}{{\Small $A$}}
  \psfrag{B}{{\Small $B$}}
  \psfrag{C}{{\Small $C$}}
  \psfrag{D}{{\Small $D$}}
  \psfrag{E}{{\Small $E$}}
  \psfrag{F}{{\Small $F$}}
  \psfrag{G}{{\Small $G$}}
  \psfrag{H}{{\Small $H$}}
  \includegraphics{twoorbits1bis.eps}
\end{center}

In this case the same vertical cylinder is to the right of $A$ and
$H$.
If the vertical separatrix going down from $H$ ends in $B$, there is
only one vertical cylinder (one horizontal cylinder for the rotated
surface $\V S$); if not, it necessarily crosses the shaded region to
the right of $B$, so there are two vertical cylinders, and the sum of
their heights is at most $t_2$ (the twist of the bottom cylinder of
$S$), hence less than the height of the bottom cylinder of $S$.

Case 3.
The bottom twist is zero but the top twist is nonzero.

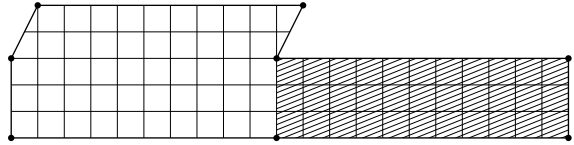
\begin{center}
\includegraphics{twoorbits2.eps}
\end{center}

Act by $\V$; this rotates $S$ by $\pi/2$.
The rotated surface $\V \cdot S$ has two cylinders: a top cylinder,
corresponding to the side part of $S$ (shaded on the figure), with
twist $0$, and a bottom cylinder cylinder of height at most $t_1$,
which we assumed to be less than $h_1$.
The surface in the same cusp with least nonnegative twists also has
top twist $0$, so if it has bottom twist $0$, conclude by case 4,
otherwise apply case 2 to obtain a surface of height less than $h_1$.

Case 4.
The twist parameters are both zero.
In this case we end the induction by jumping to a one-cylinder surface
directly:

\begin{Lemma}
The diagonal direction for the ``base rectangle'' of an L surface
tiled by a prime number of squares is a one-cylinder direction.
\end{Lemma}

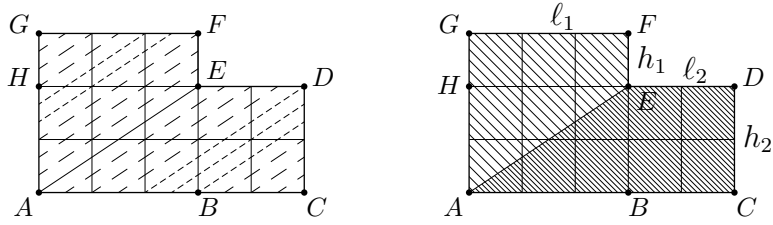
\begin{center}
  \psfrag{A}{{\Small $A$}}
  \psfrag{B}{{\Small $B$}}
  \psfrag{C}{{\Small $C$}}
  \psfrag{D}{{\Small $D$}}
  \psfrag{E}{{\Small $E$}}
  \psfrag{F}{{\Small $F$}}
  \psfrag{G}{{\Small $G$}}
  \psfrag{H}{{\Small $H$}}
  \psfrag{l}{$\ell_1$}
  \psfrag{ll}{$\ell_2$}
  \psfrag{h}{$h_1$}
  \psfrag{hh}{$h_2$}
  \includegraphics{magicdir.eps}
  \hfil
  \includegraphics{magicdirproof.eps}
\end{center}

\begin{proof}
The ascending diagonal $[AE]$ of the base rectangle of our L
surface cuts it into two zones.
Note that $[AE]$ has no other integer point than $A$ and $E$ by (P) of
\S\,\ref{param:two:cyl:surf}.

The other two saddle connections parallel to $[AE]$ start from $B$ and
$H$ and end in $F$ and $D$.
We want to prove that the one starting from $H$ ends in $F$ and the
one issued from $B$ ends in $D$, meaning each saddle connection
returns with angle $3 \pi$.

Set the origin in $A$ or $E$ and consider coordinates modulo $\ell_1
\Z \times h_2 \Z$.

Follow a saddle connection parallel to $[AE]$ from integer point to
integer point.
While it winds in a same zone, the coordinates of the integer points
it reaches remain constant modulo $\ell_1 \Z \times h_2 \Z$.
Changing zone has the following effects for the coordinates of the
next integer point:\\
\textbullet\enspace from the upper to the lower zone: decrease $y$ by
$h_1$ modulo $h_2$;\\
\textbullet\enspace from the lower to the upper zone: decrease $x$ by
$\ell_2$ modulo $\ell_1$.

Zone changes have to be alternated.
Once inside a zone with the right coordinates modulo $\ell_1 \Z \times
h_2 \Z$, a separatrix reaches the top right corner of the zone with no
more zone change.

So we want to prove that starting from $B$, in the lower zone with
coordinates $(0,0)$, and adding in turn $(- \ell_2, 0)$ and $(0, -
h_1)$, coordinates $(\ell_2, 0)$ (point $D$) will be reached before
$(0, h_1)$ (point $H$).

After $k$ changes from lower to upper zone and $k$ changes from upper
to lower zone, the coordinates are final if $k \equiv - 1
\pmod{\ell_1}$ and $k \equiv 0 \pmod{h_2}$; that is, if $k$ is
$h_2(\ell_1 - 1)$.
After $k+1$ changes from lower to upper zone and $k$ changes from
upper to lower zone, the coordinates are final if $k \equiv 0
\pmod{\ell_1}$ and $k \equiv 0 \pmod{h_2}$, which means $k$ is $h_2
\cdot \ell_1$.
So the separatrix parallel to $[AE]$ starting from $B$ reaches $D$.
\end{proof}

\subsection{Linking one-cylinder surfaces of each type}
\label{sec:linkonecyl}

We call a surface type A (resp.\ B) if it has $1$ (resp.\ $3$) integer
Weierstrass points.

Recall that a primitive one-cylinder surface in $\cH(2)$ has height
one, hence it is determined by the cyclically ordered lengths of the
three saddle connections on the bottom of this cylinder (which add up
to $n$), and by a twist parameter.

The repeated action of $\U$ can set the twist parameter to any of its
$n$ possible values, so for the purpose of linking surfaces of the
same type by $\SL(2,\Z)$ action, we may already consider surfaces with
the same cyclically ordered partition $(a,b,c)$ as equivalent
(allowing implicit $\U$-action).
We will call them $(a,b,c)$ surfaces.

Partitions into three odd numbers correspond to type A; partitions
into two even numbers and one odd number correspond to type B.

We will first show that any one-cylinder surface has a $(1,*,*)$
surface in its orbit; then we will show that $(1,b,c)$ surfaces with
$b$ and $c$ odd are in the orbit of a $(1,1,n-2)$ surface, proving all
type A surfaces to be in one orbit; then that $(1,2a,2b)$ surfaces are
in the orbit of a $(1,2,n-3)$ surface, proving all type B surfaces to
be in one orbit.

Consider a rational-slope direction on a square-tiled surface $S$;
this direction is completely periodic.
Say it is given by a vector $(p,q)\in\Z^2$, with $p\wedge q=1$.
For any $(u,v)\in\Z^2$ such that $\det\mat{p}{q}{u}{v}=1$ our surface
can be seen as tiled by parallelograms of sides $(p,q)$, $(u,v)$,
whose vertices are the vertices of the square tiling.

These parallelograms are taken to unit squares by
$M=\mat{p}{q}{u}{v}^{-1}\in\SL(2,\Z)$.
We call $M\cdot S$ ``the surface seen in direction $(p,q)$'' on $S$.

Consider a saddle connection $\sigma$ on $S$ in direction $(p,q)$; the
corresponding saddle connection on $M\cdot S$ is horizontal with an
integer length equal to the number of integer points (vertices of the
square tiling) $\sigma$ reaches on $S$.
Abusing vocabulary we also call this the length of $\sigma$.

A saddle connection returns at an angle of $3\pi$ if and only if it
has a Weierstrass point in its middle.
If two saddle connections in a given direction return with angle
$3\pi$ then so does the third, and that direction is one-cylinder;
thus two saddle connection lengths give the third.

\subsubsection{First step: any one-cylinder surface has a $(1, *, *)$
surface in its orbit}
\label{onecylfirststep}

To show this, we prove that an $(a,b,c)$ surface has a
$(\delta, k\delta, \gamma)$ surface in its orbit, where
$\delta \,|\, a \wedge b$.
Then because $n$ is prime we have $\gamma \wedge \delta = 1$, hence
applying the argument a second time with $\gamma$ and $\delta$ in
place of $a$ and $b$ shows that there is a $(1,*,*)$ one-cylinder
surface in the orbit of the surface we started with.

The proof is as follows.
Consider the $(a, b, c)$ surface $S$ having saddle connections of
lengths $a$, $b$, $c$ on the bottom, $b$, $a$, $c$ on the top.

\begin{wrapfigure}{r}{0pt}
  \psfrag{a}{$a$}
  \psfrag{b}{$b$}
  \psfrag{c}{$c$}
  \noindent
  \includegraphics{onecyl3.eps}
\end{wrapfigure}

$\V S$ has two cylinders, the top one of height $c$ and width $1$, and
the bottom one of height $d = a \wedge b$ and width $\frac{a + b}{d}$,
and some twist $t$.

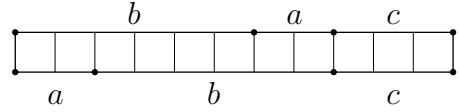
\begin{wrapfigure}{r}{0pt}
  \noindent
  \includegraphics{onecyl4.eps}
\end{wrapfigure}

Now the direction $(1 + t, d)$ is a $(\delta, k \delta, \gamma)$
one-cylinder direction with $\delta = (1 + t) \wedge d$.
Note that $k = \frac{a + b}{d} - 1$, and that $\gamma \wedge \delta =
1$.

So by applying this procedure twice we see that any surface has a $(1,
*, *)$ one-cylinder surface in its orbit.

\subsubsection{End of proof for type A surfaces}

There only remains to link any $(1, b, c)$ surface, where $b$ and $c$
are odd, to a $(1, 1, n - 2)$ surface.

\begin{wrapfigure}{r}{0pt}
  \psfrag{b}{$b$}
  \psfrag{c}{$c$}
  \psfrag{1}{$1$}
  \noindent
  \includegraphics{onecyl1.eps}
\end{wrapfigure}

Consider the L surface with arms of width $1$ and lengths $b$
and $c$.

Apply $\U^{2}$ to set the bottom twist to $2$.
Then rotate by applying $\V$, and obtain a surface with two cylinders
of height $1$.
By applying a convenient power of $\U$ the twists can be made both
$0$.

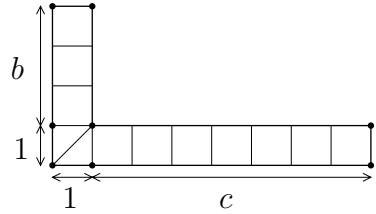
\begin{wrapfigure}{r}{0pt}
  \psfrag{b}{$b$}
  \noindent
  \includegraphics{onecyl2.eps}
\end{wrapfigure}

In the diagonal direction of the base rectangle of this new L surface,
we see a $(1, 1, n - 2)$ surface.

\subsubsection{End of proof for type B surfaces}

Here we take the one-cylinder surface with the partition $(1,2,n-3)$
as the reference surface, and prove by steps that any type B surface
has it in its orbit.

To do this, we first show that any one-cylinder surface has a
one-cylinder surface with a $(1,2a,2b)$ partition in its orbit.
This is done by the first step explained above.

Then we link\\
\textbullet\enspace $(1,2a,2b)$ where $a \neq b$ with $(d,2d,*)$, then
with $(1,2,n-3)$;\\
\textbullet\enspace $(1,2a,2b)$ where $a = b$ with $(2,2,n-4)$, then
with $(1,2,n-3)$.

\begin{itemize}
\item Linking  $(1,2a,2b)$ with $(1,2,*)$ when $a \neq b$.
\end{itemize}

Without loss of generality, suppose $a < b$.
Consider the one-cylinder surface with saddle connections of lengths
$2a$, $2b$, $1$ on the bottom and $2b$, $2a$, $1$ on the top.

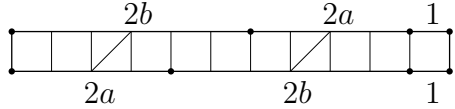
\begin{center}
  \psfrag{1}{1}
  \psfrag{2a}{$2a$}
  \psfrag{2b}{$2b$}
  \includegraphics{onecyl5.eps}
\end{center}

In the direction $(b-a, 1)$ there is a connection between two integer
Weierstrass points, so in this direction we see a two-cylinder
surface.
Its top cylinder has height $2a$ and width $2$ and its bottom cylinder
has height $1$ and with $2 + \ell$ for some $\ell$.

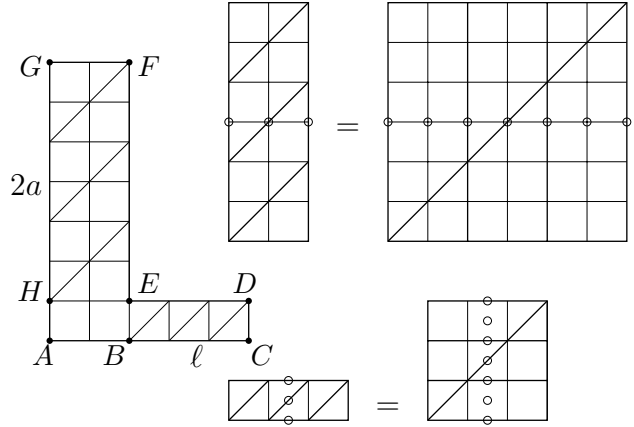
\begin{center}
  \psfrag{A}{{\small $A$}}
  \psfrag{B}{{\small $B$}}
  \psfrag{C}{{\small $C$}}
  \psfrag{D}{{\small $D$}}
  \psfrag{E}{{\small $E$}}
  \psfrag{F}{{\small $F$}}
  \psfrag{G}{{\small $G$}}
  \psfrag{H}{{\small $H$}}
  \psfrag{l}{$\ell$}
  \psfrag{2a}{$2a$}
  \psfrag{e}{$=$}
  \includegraphics{onecyl6.eps}
\end{center}

In certain directions, the separatrix issued from $H$ winds around the
horizontal cylinder $HEGF$.
In particular, in any direction $(k, a)$, $k\in\N$, it will run into a
Weierstrass point (and into a saddle after twice the distance).

Likewise, in appropriate directions, the separatrix issued from $B$
winds around the vertical cylinder $BCDE$.
In particular, in any direction $(\ell/2, k/2)$ (equivalently $\ell,
k$), $k\in\N$, it will run into a Weierstrass point (and into a saddle
after twice the distance).

Consider therefore the direction $(\ell, a)$.
In this direction we get a $(d, 2d, *)$ one-cylinder surface, where $d
= a \wedge \ell$.

Now there only remains to link $(d, 2d, *)$ with $(1, 2, *)$, which is
easily done: consider the one-cylinder surface with saddle connections
$d$, $2d$, $c$ on the bottom and $2d$, $d$, $c$ on the top;

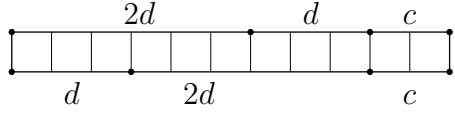
\begin{center}
  \psfrag{d}{$d$}
  \psfrag{2d}{$2d$}
  \psfrag{c}{$c$}
  \includegraphics{onecyl7.eps}
\end{center}

\noindent in the $(d, 1)$ direction we get a $(1, 2, *)$ one-cylinder
surface.

\begin{itemize}
\item Linking  $(1,2a,2b)$ with $(1,2,*)$ when $a = b$.
\end{itemize}

Consider the one-cylinder surface with saddle connections of length
$2a$, $2b$, $c$ on the bottom and $2b$, $2a$, $c$ on the top.

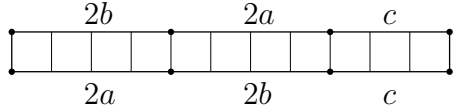
\begin{center}
  \psfrag{c}{$c$}
  \psfrag{2a}{$2a$}
  \psfrag{2b}{$2b$}
  \includegraphics{onecyl8.eps}
\end{center}

In the direction $(a,1)$ we see a $(2,2,*)$ one-cylinder surface.

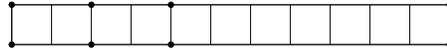
\begin{center}
\includegraphics{onecyl9.eps}
\end{center}

On this surface, in the direction $(2,1)$, we have a two-cylinder
surface with its top cylinder of height $2$ and width $1$, and its
bottom cylinder of height $1$.
Acting by $\U$ we can set the twist parameters to $0$.

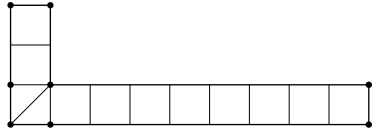
\begin{center}
\includegraphics{onecyl10.eps}
\end{center}

Then in the direction $(1, 1)$ we see a $(1, 2, n - 3)$ one-cylinder
surface.

\subsection{L-shaped billiards}

L-shaped billiards give rise to L-shaped translation surfaces by an
unfolding process; any L-shaped translation (with zero twists) surface
is the covering translation surface of an L-shaped billiard.

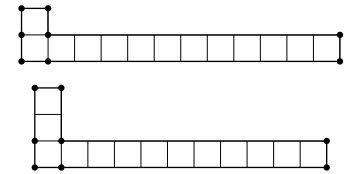
\begin{wrapfigure}{l}{0pt}
  \noindent
  \includegraphics{Lbilliardsv.eps}
\end{wrapfigure}

Fix some prime $n>3$, and consider the two-cylinder surfaces $S_1$ and
$S_2$, both having $h_2=1$, $w_1=1$ and $t_1=t_2=0$, and $S_1$ having
$h_1=1$, $w_2=n-1$ and $S_2$ having $h_2=2$, $w_2=n-2$.
The picture on the side represents $S_1$ and $S_2$ for $n = 13$.

For each $n$, $S_1$ and $S_2$ belong to orbit A and B respectively,
and arise from L-shaped billiards.
This proves Proposition~\ref{prn:L}.

\section{Proof of results about elliptic points}

Some constructions in this section are inspired by \cite{Ve95}.

\subsection{Translations}

Here we prove Proposition~\ref{prn:no:translations}.

Suppose a surface $S\in \cH(2)$ has a nontrivial translation $f$ in its
affine group.
$f$ fixes the saddle and induces a permutation on outgoing horizontal
separatrices.
Let $\varepsilon$ be smaller than the length of the shortest saddle
connection of $S$, and consider the three points at distance
$\varepsilon$ from the saddle on the three separatrices in a given
direction.
$f$ cannot fix any of these points, otherwise it would be the identity
of $S$, but it fixes the set of these points, hence it induces a
cyclic permutation on them.
This implies that except for the saddle, which is fixed, all
$f$-orbits have size $3$.
However the set of regular Weierstrass points is also fixed (since the
translation $f$ is an automorphism of the underlying Riemann surface),
and has size $5$.
This is a contradiction.

\subsection{Elliptic points of order $3$}

Here we prove Proposition~\ref{prn:ellipticthree}.

Suppose a surface $S$ in $\cH(2)$ has an elliptic element of projective
order three in its Veech group.
Since the hyperelliptic involution has order $2$, $S$ has in fact an
elliptic element of order $6$ in its Veech group.
Conjugate by $\SL(2,\R)$ to a surface that has the rotation by $\pi/3$
(hereafter denoted by $r$) in its Veech group.

Considering Proposition~\ref{prn:no:translations}, we denote by $r$
the corresponding affine diffeomorphism.

The set of Weierstrass points is preserved by $r$.
The saddle being fixed, the remaining five Weierstrass points are
setwise fixed, so at least two of them are also fixed.
Consider one Weierstrass point that is fixed, call it $W$.
Consider the shortest saddle connections through $W$.
They come by triples making angles $\pi/3$.

Take one such triple, consider the corresponding regular hexagon
(which has these saddle connections as its diagonals).

We can take this hexagon as a building block for a polygonal
fundamental domain of the surface.
Consider a pair of opposite sides of this hexagon; they cannot be
identified, since the rotational symmetry would imply other
identifications and mean we have a torus.

Hence, these sides and the diagonal parallel to them are three saddle
connections in the same direction.
So this is a completely periodic direction, and we want to see two
cylinders in this direction.
This would imply identifying two opposite sides, which we have
excluded.

\subsection{Elliptic elements of order $2$}

\subsubsection{Proof of Lemma~\ref{lemma:octagon}}

Here, inspired by \cite{Ve95}, we give a convenient representation for
$\V$-invariant Veech surfaces in $\cH(2)$: a fundamental octagon which
is $\V$-invariant.
Consider a Veech surface in $\cH(2)$ that has $\V$ in its Veech group;
denote also by $\V$ the corresponding affine diffeomorphism.

The set of Weierstrass points is fixed by $\V$ (as by any affine
diffeomorphism).
The saddle being fixed, at least one of the remaining 5 Weierstrass
points must be fixed.

Consider such a point and the shortest saddle connections through this
point.
They come by orthogonal pairs.
Take one such pair.
Consider the square having this pair of saddle connections as
diagonals.
Without loss of generality, consider the sides of the square as
horizontal and vertical.

This square is the central piece of our fundamental domain.
Other than the corners (the saddle) and the center, there are no
Weierstrass point inside this square or on its edges.

Consider the horizontal sides of our square.
These sides are saddle connections so they define a completely
periodic direction on the surface.

These sides are not identified, otherwise by $\V$-symmetry the other
two would also be and we would have a torus.
So this is a two-cylinder direction and our two sides bound the short
cylinder in this direction.
This short cylinder lies outside the square and can be represented as
a parallelogram with its ``top-left'' corner in the vertical strip
defined by the square (i.e.\ with a ``reasonable'' twist).

By $\V$-symmetry there also is such a parallelogram in the other
direction.
To make the picture more symmetric each parallelogram can be cut into
two triangles, glued to opposite sides of the square.
Thus we get a representation of the surface as an octagon with
(parallel) opposite sides identified.
Note that the four remaining Weierstrass points are the middle of the
sides of this octagon.

\subsubsection{Proof of Proposition~\ref{prn:elliptictwo}}
\label{square:octagon}

Represent the surface as above: an octagon made of a square and four
triangles glued to its sides.
All vertices lie on integer points.

\begin{wrapfigure}{r}{0pt}
\noindent
  \psfrag{A}{$A$}
  \psfrag{B}{$B$}
  \psfrag{C}{$C$}
  \includegraphics{elliptic4.eps}
\end{wrapfigure}

Let $ABC$ be one of the triangles, labeled clockwise so that $AC$ is a
side of the square.

Let $(p, q)$ be the coordinates of $\overrightarrow{AC}$ and $(r, s)$
those of $\overrightarrow{AB}$.
The area of the surface is then $p ^ 2 + q ^ 2 + 2 (p s - q r)$.

If $n$ is prime then $p$ and $q$ have to be relatively prime, and of
different parity.
Then $p ^ 2 + q ^ 2 \equiv 1 \pmod 4$.
The center of the square lies at the center of a square of the tiling.
The condition for two Weierstrass points to lie on integer points is
for $(p s - r q)$ to be even.

We conclude by observing that $n$ is $1$ (resp.\ $3$) modulo $4$ when
$(p s - r q)$ is even (resp.\ odd).

\section{Proof of countings}

Here we establish the countings and estimates of
Proposition~\ref{prn:countings}.

\subsection{One-cylinder cusps and surfaces}

For prime $n > 3$, one-cylinder $n$-square-tiled cusps in $\cH(2)$ are
in 1-1 correspondence with cyclically ordered $3$-partitions of $n$.

Ordered $3$-partitions $(a,b,c)$ of $n$ are in 1-1 correspondence with
pairs of distinct integers $\{\alpha,\beta\}$ in $\{1, \ldots, n-1\}$:
assuming $\alpha < \beta$, the correspondence is given by $a =
\alpha$, $a + b = \beta$, $a + b + c = n$.
So there are $C_{n-1}^2$ ordered $3$-partitions of $n$.
Ordered $3$-partitions of $n$ being in 3-1 correspondence with
cyclically ordered $3$-partitions, there are $\frac{1}{3}C_{n-1}^2 =
\frac{(n-1)(n-2)}{6}$ cyclically ordered $3$-partitions of $n$.

Thus there are $\frac{(n-1)(n-2)}{6}$ one-cylinder cusps of
$n$-square-tiled translation surfaces in $\cH(2)$.

Those in orbit A are those with $3$ odd parts $2a-1$, $2b - 1$, $2 c -
1$; these are in 1-1 correspondence with cyclically ordered partitions
$a$, $b$, $c$, of $\frac{n+3}{2}$.
Their number is hence $\frac{1}{6}(\frac{n+3}{2}-1)(\frac{n+3}{2}-2) =
\frac{(n+1)(n-1)}{24}$.

The remaining ones are in orbit B, their count is hence the
difference, $\frac{(n-1)(n-3)}{8}$.

All one-cylinder cusps discussed here have width $n$ ($n$ possible
values of the twist parameter), so the counts of one-cylinder surfaces
are $n$ times the corresponding cusp counts.

\subsection{Two-cylinder surfaces}
\label{sec:count:two:cyl:surf}
The total number of two-cylinder $n$-square-tiled surfaces ($n$ prime)
is
$$
S(n) = \sum_{a, b, k, \ell} k \ell,
$$
where the sum is over $a,b,k,\ell \in \N^*$ such that
$k < \ell$ and $a k + b \ell = n$.

This follows from the parametrization in \S\,\ref{param:two:cyl:surf};
the letters $a$, $b$, $k$, $\ell$ used here correspond to the
parameters $h_1$, $h_2$, $w_1$, $w_2$ there, and the summand is the
number of possible values of the twist parameters, given the heights
and widths of the two cylinders.

We want the asymptotic for this quantity as $n$ tends to infinity, $n$
prime.
In order to find this, we consider the sum as a double sum: the sum
over $a$ and $b$ of the sum over $k$ and $\ell$.

Write $S(n) = \sum_{a, b} S_{a, b}(n)$, where $S_{a, b}(n) = \sum_{k,
\ell} k \ell$.

\bigskip

We study the inner sum by analogy with a payment problem: how many
ways are there to pay $n$ units with coins worth $a$ and $b$ units?

This problem is classically solved by the use of generating series:
denote the number of ways to pay by $s_{a, b}(n)$; then

\centerline{$s_{a, b}(n) =
\Card \{ ( k, \ell ) \in \N ^ 2 : a k + b \ell = n \} = \sum_{k,
\ell \in \N : a k + b \ell = n} 1$.}

Now notice that $\sum_{k = 0}^{\infty} z^{a k} \sum_{\ell =
0}^{\infty} z^{b \ell} = \sum_{n = 0}^{\infty} s_{a, b}(n) z^n$, and
deduce that the number looked for is the $n$-th coefficient of the
power series expansion of the function $\frac{1}{1 - z^a} \frac{1}{1 -
z^b}$.

\bigskip

We turn back to our real problem, $S_{a, b}(n) = \sum_{k, \ell \in
\N^{*} : a k + b \ell = n, k < \ell} k \ell$.

We want to show that $S(n) \sim c n^3$ for prime $n$.
For this we will use the dominated convergence theorem: we show that
$S_{a, b}(n) / n^3$ has a limit $c_{a, b}$ when $n$ tends to infinity
with $a$ and $b$ fixed, and that $S_{a, b}(n)/n^3$ is bounded by some
$g_{a, b}$ such that $\sum_{a,b}g_{a, b}<\infty$, to conclude that
$S(n)/n^3 = \sum_{a, b}S_{a, b}(n)/n^3$ tends to $c = \sum_{a, b}c_{a,
b}$, which means $S(n) \sim c n^3$.

The dominated convergence is proved as follows.

Write $S_{a, b}(n) = \sum_{k, h \in \N^{*} : (a + b) k + b h = n} k (k
+ h)$ by introducing $h = \ell - k$.
Then split the sum into $\sum k^2$ and $\sum kh$.
Write
$$
S'_{a, b}(n) = \sum_{k,h\in\N^*,\ (a + b) k + b h = n} k^2/{n^3}
\leqslant \sum_{k\in\N^*,\ h\in\Q,\ (a + b) k + b h = n} k^2/{n^3}
$$
$$
S''_{a, b}(n) =
\sum_{k,h\in\N^*,\ (a + b) k + b h = n} kh/{n^3}
\leqslant
\sum_{k\in\N^*,\ h\in\Q,\ (a + b) k + b h = n} kh/{n^3}
$$
(in the sums on the right-hand side, $h$ has been allowed to be a
rational instead of an integer.)
Hence
$$
S'_{a, b}(n) \leqslant \frac{1}{(a+b)^3} \Bigl[\frac{a+b}{n}
\sum_{k=1}^{\lfloor n/(a+b)\rfloor} \bigl(\frac{a+b}{n}
k\bigr)^2\Bigr]
$$
$$
S''_{a, b}(n) \leqslant
\frac{1}{(a+b)^2 b} \Bigl[\frac{a+b}{n} \sum_{k=0}^{n/(a+b)}
\bigl(\frac{a+b}{n} k\bigr) \bigl(1 - \frac{a+b}{n} k\bigr)\Bigr]
$$
The expressions in brackets, Riemann sum approximations to the
integrals $\int_{0}^{1}x^2 \d x$ and $\int_{0}^{1}x (1 - x) \d x$, are
uniformly bounded by $1$.

Now notice that $\sum_{a,b} \frac{1}{(a+b)^3}$ and $\sum_{a,b}
\frac{1}{(a+b)^2 b}$ are convergent.
This ends the dominated convergence argument.

We can now investigate the limit.
For ease of calculation, we drop the condition $k< \ell$.
We take care of it by writing $\sum_{k,\ell} = 2 \sum_{k<\ell} +
\sum_{k = \ell}$.
For prime $n$, $k = \ell$ implies that they are both equal to $1$.
The sum for $k = \ell$ is hence equal to $n - 1$, and we will not need
to take it into account since the whole sum will grow as $n^3$.

Denote by $\tilde S(n,a,b)$ the sum over all $k$ and $\ell$.

Notice that $\sum_{k = 0}^{\infty} k z^{a k} \sum_{\ell = 0}^{\infty}
\ell z^{b \ell} = \sum_{n = 0}^{\infty} \tilde S(n, a, b) z^n$.

$\tilde S(n,a,b)$ is therefore the $n$-th coefficient of the power
series expansion of the function $f_{a,b} =
\frac{z^a}{(1-z^a)^2}\frac{z^b}{(1-z^b)^2}$.

To determine this coefficient, decompose $f_{a,b}$ into partial
fractions.
This function has poles at $a$-th and $b$-th roots of $1$.
Since $n$ is prime, we are only interested in relatively prime $a$ and
$b$, for which the only common root of $1$ is $1$ itself, which is
hence a $4$-th order pole of $f_{a,b}$, while other poles have order
$2$.

The $n$-th coefficient of the power series expansion of $f_{a,b}$ is a
polynomial of degree $3$ in $n$, whose leading term is $c_{a,b}
\frac{n^3}{6}$, where $c_{a,b}$ is the coefficient of
$\frac{1}{(1-z)^4}$ in the decomposition of $f_{a,b}$ into partial
fractions.
This coefficient is computed to be $\frac{1}{a^2b^2}$.

We want the sum over relatively prime $a$ and $b$.
We relate it to the sum over all $a$ and $b$ by sorting the latter
according to $d = a \wedge b$.
$$
\sum_{a,b}\frac{1}{a^2b^2} = \sum_{d}\sum_{a,b,\ a \wedge b =
d}\frac{1}{a^2b^2} = \sum_{d}\frac{1}{d^4}\sum_{a,b,\ a \wedge b =
1}\frac{1}{a^2b^2}.
$$

By observing that $\sum_{a,b}\frac{1}{a^2b^2} =
(\sum_{a}\frac{1}{a^2})^2 = \zeta(2)^2 = \frac{\pi^4}{36}$ and that
$\sum_{d}\frac{1}{d^4} = \zeta(4) = \frac{\pi^4}{90}$ we get that the
sum $\sum_{a,b,\ a \wedge b = 1}\frac{1}{a^2b^2}$ is equal to $5/2$.
Divide by $2$ to get back to $k<\ell$, and find that $S(n) \sim
\frac{5}{4}\frac{n^3}{6}$.

\subsection{Two-cylinder surfaces by orbit}

Two-cylinder surfaces for which both heights are odd are in orbit A;
those for which both widths are odd are in orbit B; half of the
remaining ones are in orbit A, and half in B; the factor one half
comes from the conditions on the twists.
(See the table in \S\,\ref{sec:twovalinv}.)

First compute the asymptotic for \textbf{odd heights}.
Write
$$
S^\mathrm{oh}(n) = \sum_{\substack{a,b,k,\ell \\ ak + b\ell = n \\ a,
b\text{ odd} \\ a \wedge b = 1 \\ k < \ell}}k\ell.
$$
Then $S^\mathrm{oh}(n) \sim \frac{1}{2} \widetilde S^\mathrm{oh}(n)$
where $\widetilde S^\mathrm{oh}(n)$ is the same sum without the
condition $k < \ell$.
The dominated convergence works as previously.

For odd $a$ and $b$ such that $a \wedge b = 1$,
$$
\widetilde S^\mathrm{oh}_{a,b}(n) = \sum_{\substack{k,\ell \\ a k + b
\ell = n}} k \ell \sim \frac{1}{a^2b^2}\cdot\frac{n^3}{6}.
$$

We need to sum over relatively prime odd $a$ and $b$.
Using the same trick as previously, write
$$
\sum_{a, b \text{ odd}}\frac{1}{a^2b^2} = \sum_{d \text{ odd}}
\sum_{\substack{a, b \text{ odd} \\ a\wedge b = d}} \frac{1}{a^2b^2} =
\sum_{d \text{ odd}}\frac{1}{d^4} \sum_{\substack{a, b \text{ odd} \\
a\wedge b = 1}} \frac{1}{a^2b^2}.
$$
Now
$$
\sum_{a, b \text{ odd}}\frac{1}{a^2b^2} = \Bigl(\sum_{a \text{
odd}}\frac{1}{a^2}\Bigr)^2 = ((1 - 1/2^2)\zeta(2))^2 =
\frac{9}{16}\cdot\frac{\pi^4}{36}
$$
and
$$
\sum_{d \text{ odd}}\frac{1}{d^4} = (1 - 1/2^4)\zeta(4) =
\frac{15}{16}\cdot\frac{\pi^4}{90}
$$
so
$$
\sum_{\substack{a, b \text{ odd} \\ a\wedge b = 1}}
\frac{1}{a^2b^2} = 3/2.
$$

We deduce that $S^\mathrm{oh}(n) \sim \frac{3}{4}\,\frac{n^3}{6}$ (the
condition $k < \ell$ is responsible for a factor $1/2$).

Similarly compute the asymptotic for \textbf{odd widths}.
Write
$$
S^\mathrm{ow}(n) = \sum_{\substack{a,b,k,\ell \\ ak + b\ell = n \\ k,
\ell\text{ odd} \\ a \wedge b = 1 \\ k < \ell}}k\ell.
$$
For fixed $a$ and $b$ with $a \wedge b = 1$, put
$$
\widetilde S^\mathrm{ow}_{a,b}(n) = \sum_{\substack{k,\ell \text{
odd}\\ a k + b \ell = n}} k \ell.
$$

Notice that $\sum_{k\text{ odd}}kz^{ak}\sum_{\ell\text{ odd}}\ell
z^{b\ell} = \sum \widetilde S^\mathrm{ow}_{a,b}(n)z^n$.

Because $\sum kz^k = \frac{z}{(1-z)^2}$, $\sum 2kz^{2k} =
\frac{2z^2}{(1-z^2)^2}$, and the difference is $\sum (2k + 1)z^{2k+1}
= \frac{z(1+z^2)}{(1-z^2)^2}$.

$\widetilde S^\mathrm{ow}_{a,b}(n)$ is now the $n$-th coefficient of
the power series expansion of $\frac{z^a(1+z^{2a})}{(1-z^{2a})^2}
\cdot \frac{z^b(1+z^{2b})}{(1-z^{2b})^2}$.
When $a\wedge b = 1$, this rational function has two order $4$ poles
at $1$ and $-1$ and its other poles have order $2$; the coefficients
of $\frac{1}{(1-z)^4}$ and $\frac{1}{(1+z)^4}$ in its decomposition
into partial fractions are respectively $\frac{1}{4a^2b^2}$ and
$\frac{(-1)^{a+b}}{4a^2b^2}$.

Because $n$ is odd, and $k$ and $\ell$ are odd, $a$ and $b$ have to
have different parities, so $a + b$ is odd.
So $\frac{1}{4a^2b^2} + \frac{(-1)^{a+b}(-1)^n}{4a^2b^2} =
\frac{1}{2a^2b^2}$.

Now
$$
\sum_{\substack{a \wedge b = 1 \\ a \not \equiv b
\pmod{2}}}\frac{1}{a^2b^2} = \sum_{a \wedge b = 1}\frac{1}{a^2b^2} -
\sum_{\substack{a \wedge b = 1 \\ a, b \text{ odd} }}\frac{1}{a^2b^2}
= 5/2 - 3/2 = 1.
$$

The condition $k < \ell$ brings a factor $1/2$, thus we get
$S^\mathrm{ow}(n)\sim (1/4)(n^3/6)$.

The remaining surfaces are those for which heights as well as widths
of the cylinders have mixed parities.
The asymptotic for this ``\textbf{even-odd}'' part is computed as the
difference between the total sum and the odd-widths and odd-heights
sums.

Write $S(n) = S^\mathrm{oh}(n) + S^\mathrm{ow}(n) + S^\mathrm{eo}(n)$.
We already know that $S(n) \sim \frac{5}{4}\cdot\frac{n^3}{6}$,
$S^\mathrm{oh}(n) \sim \frac{3}{4}\cdot\frac{n^3}{6}$, and
$S^\mathrm{ow}(n) \sim \frac{1}{4}\cdot\frac{n^3}{6}$.
So the even-odd part has asymptotics $S^\mathrm{eo}(n) \sim
\frac{1}{4}\cdot\frac{n^3}{6}$.

Putting pieces together, the number of $n$-square-tiled two-cylinder
surfaces of type A, $n$ prime, is equivalent to $(3/4 + 1/8)(n^3/6) =
(7/8)(n^3/6)$.
For type B, we get $(1/4 + 1/8)(n^3/6) = (3/8)(n^3/6)$.

\subsection{Two-cylinder cusps}

For $n$ prime, the number of two-cylinder cusps (in both orbits) is
given by
$$
S(n) = \sum_{\substack{a, b, k, \ell \in \N^*\\ ak+b\ell=n \\
k<\ell}} (a \wedge k) (b \wedge \ell).
$$
(see counting of two-cylinder surfaces in
\S\,\ref{sec:count:two:cyl:surf} and discussion of cusps in
\S\,\ref{sec:tools:action}.)

\begin{Remark}
For nonprime $n$, the number of two-cyl cusps is less than $S(n)$
defined as above, so the bound found here is still valid.
\end{Remark}

$S(n)$ is less than
$$
\widetilde S(n) = \sum_{\substack{a, b, k, \ell \in \N^*\\
ak+b\ell=n}} (a \wedge k) (b \wedge \ell).
$$
where the condition $k < \ell$ is dropped.

We will show that for any $\varepsilon>0$, $\widetilde S(n)
\ll_{n\to\infty} n ^{3/2 + \varepsilon}$.

This will imply that the number of two-cylinder cusps of
$n$-square-tiled surfaces is sub-quadratic, thus negligible before the
(quadratic) number of one-cylinder cusps in each orbit.

$$
\widetilde S(n) = \sum_{\substack{A,B,u,v \in \N^*\\ Au^2+Bv^2=n}} uv
f(A) f(B),\quad\text{where }f(m) = \displaystyle\sum_{\substack{rs =
m\\ r \wedge s=1}}1.
$$

Note that $f(m) \leqslant \d(m) \ll m^\varepsilon$, where $\d(m)$ is
the number of divisors of $m$.
The factors $f(A)f(B)$ therefore contribute less than an
$n^\varepsilon$.

$$
\sum_{\substack{A,B,u,v \in \N^*\\ Au^2+Bv^2=n}} uv = \sum_{u} u
\sum_{A\leqslant n/ u^2}\Bigl(\sum_{v^2 | n-A u^2}v \Bigr).
$$

The sum in parentheses has less than $\d(n-Au^2)$ summands, each of
which is bounded by $\sqrt{n-Au^2}$, so
$$
\widetilde S(n) \ll n^{1/2+2\varepsilon} \sum_{u} n/u \ll
n^{3/2+3\varepsilon}.
$$
\hfill \qed

We thank Jo\"el Rivat for contributing this estimate \cite{Ri}.

\subsection{Elliptic points}

The discussion in \S\,\ref{square:octagon} implies that their
number is less than the number of integer-coordinate vectors in a
quarter of a circle of radius $\sqrt n$, so it is $O(n)$.

\section{Strong numerical evidence}
\label{section:numerical:evidence}

Martin Schmoll pointed out to us that the number of primitive
$n$-square-tiled surfaces in $\cH(2)$ is given in \cite{EsMaSc} to be
$$
\frac{3}{8}(n-2)n^2\prod_{p|n}(1-\frac{1}{p^2}).
$$
By \cite{Mc2}, for even $n$ all these surfaces are in the same orbit,
and for odd $n \geqslant 5$ they fall into two orbits.
So Eskin, Masur and Schmoll's formula gives the cardinality of the
single orbit for even $n$, and the sum of the cardinalities of the two
orbits for odd $n$.

\begin{Conjecture}
\label{conj:countings}
For odd $n$, the cardinalities of the orbits are given by the
following functions:\\
orbit A: $\frac{3}{16}(n-1)n^2\prod_{p|n}(1-\frac{1}{p^2})$,\\
orbit B:  $\frac{3}{16}(n-3)n^2\prod_{p|n}(1-\frac{1}{p^2})$.
\end{Conjecture}

These formulae give degree $3$ polynomials when restricted to prime
$n$, for which Theorem~\ref{thm:countings} gives the leading term.
These polynomials are expressed in the table below.

\smallskip

\begin{small}
\noindent
\begin{tabular}{@{}lccc@{}}
  \hline
  & one-cylinder & two-cylinder & all \\
  \hline
  A
  & $\frac{1}{24}(n ^ 3 - n)$
  & $\frac{1}{48}(7 n ^ 3 - 9 n ^ 2 - 7 n + 9)$
  & $\frac{3}{16}(n ^ 3 - n ^ 2 - n + 1)$ \\[1pt]
  \hline
  B
  & $\frac{1}{8}(n ^ 3 - 4 n ^ 2 + 3 n)$
  & $\frac{1}{16}(n ^ 3 - n ^ 2 - 9 n + 9)$
  & $\frac{3}{16}(n ^ 3 - 3 n ^ 2 - n + 3)$ \\[1pt]
  \hline
  all
  & $\frac{1}{6}(n ^ 3 + 3 n ^ 2 + 2 n)$
  & $\frac{1}{24}(5 n ^ 3 - 6 n ^ 2 - 17 n + 18)$
  & $\frac{3}{8}(n ^ 3 - 2 n ^ 2 - n + 2)$ \\[1pt]
  \hline
\end{tabular}
\end{small}

\smallskip

On the other hand, the counting functions for two-cylinder cusps are
not polynomials.

\begin{Conjecture}
For prime $n$, the number of elliptic points is $\lfloor \frac{n +
1}{4}\rfloor$.
\end{Conjecture}

This conjecture is valid for the first thousand odd primes.

\appendix

\section{$n=3$ and $n=5$}
\label{app:n:3:5}

\subsection*{$n=3$}

For $n = 3$, we have the following three surfaces.

\begin{wrapfigure}{l}{0pt}
  \noindent
  \includegraphics{surfaces3.eps}
\end{wrapfigure}

If we call $S_1$ the one-cylinder surface, and $S_2$ and $S_3$ the
two-cylinder surfaces, the generators of $\SL(2,\Z)$ act as follows:
$\U S_1 = S_1$, $\U S_2 = S_3$, $\U S_3 = S_2$, $\V S_1 = S_3$, $\V
S_2 = S_2$, $\V S_3 = S_1$.
So there is only one orbit, containing $d = 3$ surfaces, the number of
cusps is $c = 2$, the number of elliptic points ($\V$-invariant
surfaces) is $e = 1$, so the genus is $g = 0$ by the Gauss-Bonnet
formula.

\subsection*{$n=5$}

For $n = 5$, we have $27$ surfaces forming $8$ cusps, a representative of
which appears on the following picture.

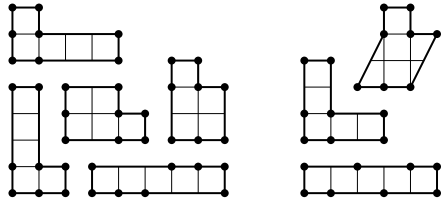
\begin{center}
  \includegraphics{cusps5.eps}
\end{center}

Computing the $\SL(2,\Z)$ action shows that they fall into two orbits,
orbit A being made of the surfaces on the left and orbit B of those on
the right.

The data for orbit A is $d = 18$ surfaces, $c = 5$ cusps, $e = 0$
elliptic point, so the genus is $g = 0$ by the Gauss-Bonnet formula.

The data for orbit B is $d = 9$ surfaces, $c = 3$ cusps, $e = 1$
elliptic point, so the genus is $g = 0$ by the Gauss-Bonnet formula.

By inspection of the congruence subgroups of genus $0$ of $\SL(2,\Z)$
(see for example \cite{CuPa}), the stabilizers of orbits A and B are
noncongruence subgroups of $\SL(2,\Z)$.

\section{Hyperelliptic components of other strata}
\label{app:hypel}

For all hyperelliptic square-tiled surfaces, one can count the number
of Weierstrass points with integer coordinates.
This provides an invariant for the action of $\SL(2,\Z)$ on
square-tiled surfaces in all hyperelliptic components of strata of
moduli spaces of abelian differentials.

The strata with hyperelliptic components are $\cH(2 g - 2)$ and $\cH(g
- 1, g - 1)$, for $g > 1$.

\begin{Proposition}
In $\cH(2 g - 2) ^ \mathrm{hyp}$ and $\cH(g - 1, g - 1) ^
\mathrm{hyp}$, for large enough odd $n$ there are at least $g$ orbits
containing one-cylinder surfaces.
\end{Proposition}

This is proved by the following reasoning.

Completely periodic surfaces in $\cH(2 g - 2)$ or $\cH(g - 1, g - 1)$,
for $g > 1$, have respectively $2 g - 1$ and $2 g$ saddle connections.

For one-cylinder primitive surfaces (necessarily of height $1$), the
lengths of the saddle connections add up to $n$, and the Weierstrass
points are two points on the circle at half-height of this cylinder
(these do not have integer coordinates), the saddle in the $\cH(2 g -
2) ^ \mathrm{hyp}$ case, and the midpoints of the saddle connections
that bound the cylinder (these have integer coordinates for exactly
those saddle connections of even length).

If $n$ is odd, the sum of the lengths is odd.
So the number of odd-length saddle connections has to be odd, and is
between $1$ and $2 g - 1$.
There are $g$ possibilities for that.
Since the value of the invariant is the number of even-length saddle
connections, it can take $g$ different values.

\section{The theorem of Gutkin and Judge}
\label{app:GuJu}

\begin{Theorem*}[Gutkin--Judge]
$(S,\omega)$ has an arithmetic Veech group if and only if $(S,\omega)$
is parallelogram-tiled.
\end{Theorem*}

Up to conjugating by an element of $\SL(2,\R)$, it suffices to show:

\begin{Theorem*}
$(S,\omega)$ is a square-tiled surface if and only if $V(S,\omega)$ is
commensurable to $\SL(2,\Z)$.
\end{Theorem*}

\noindent (i.e.\ these two groups share a common subgroup of
finite index in each.)

\begin{Remark}
In this theorem, the size of the square tiles is not assumed to be
$1$.
One can always act by a homothety to make this true, and we will
suppose that in the proof of the direct way of this theorem.
\end{Remark}

\subsection{A square-tiled surface has an arithmetic Veech group}
\label{sec:sts:arith}

Consider a square-tiled surface $(S,\omega)$, and its lattice of
periods $\Lambda(\omega)$.
By Lemma~\ref{Schmoll:Zorich}, $V(S,\omega) <
V(\rquotient{\R^2}{\Lambda(\omega)},\d z)$.

Case 1.
Let us first assume that $\Lambda(\omega)=\Z^2$, i.e.\ $(S,\omega)$ is
a primitive square-tiled surface.

Lemma~\ref{primitive:orbit} implies that $\SL(2,\Z)$ acts on the set
$E$ of square-tiled surfaces contained in its $\SL(2,\R)$-orbit.
The set $E$ is finite and the stabilizer of this action is
$V(S,\omega)$.
The class formula then implies that $V(S,\omega)$ has finite index in
$\SL(2,\Z)$.

Case 2.
Suppose that $\Lambda(\omega)$ is a strict sublattice of $\Z^2$.
Consider $P_1, \ldots, P_k$ the preimages of the origin on $S$.
Denote by $\Aff_{P_1,\ldots,P_k}$ the stabilizer of the set of these
points in the affine group of $(S,\omega)$, and $V(P_1,\ldots,P_k)$
the associated Veech group.
The translation surface $(S,\omega,\{P_1,\ldots,P_k\})$ where
$\{P_1,\ldots,P_k\}$ are artificially marked is a primitive
square-tiled surface.
From Case~1 above, its Veech group $V(P_1,\ldots,P_k)$ is therefore a
lattice contained in the discrete group $V(S,\omega)$, hence of finite
index in this group.

Thus $V(P_1,\ldots,P_k)$ is a finite-index subgroup in both
$V(S,\omega)$ and $\SL(2,\Z)$.

\subsection{A surface with an arithmetic Veech group is square-tiled}

This part is inspired by ideas of Thurston \cite{Th} and Veech
\cite[\S 9]{Ve89}, and appeared in \cite[appendix B]{Hu}.

Let $S$ be a translation surface with an arithmetic Veech group
$\Gamma$.

If $\Gamma$ is commensurable to $\SL(2,\Z)$ only in the wide sense, we
move to the case of strict commensurability.
This conjugacy on Veech groups is obtained by $\SL(2,\R)$ action on
surfaces.

We prove the following propositions.

\begin{Proposition}
\label{prn:arith:cont:parab}
A group $\Gamma$ commensurable with $\SL(2,\Z)$ contains two elements
of the form $\mat{1}{0}{m}{1}$ and $\mat{1}{n}{0}{1}$ for some
$m,n\in\N^{*}$.
\end{Proposition}

\begin{Proposition}
\label{prn:cont:parab:squ:tiled}
If the Veech group $\Gamma$ of a translation surface $S$ contains two
elements of the form $\mat{1}{0}{m}{1}$ and $\mat{1}{n}{0}{1}$ for
some $m,n\in\N^{*}$, then $S$ is square-tiled.
\end{Proposition}

Proposition~\ref{prn:arith:cont:parab} follows from the following
lemma.

\begin{Lemma}
If $H\leq G$ is a finite-index subgroup then every $g\in G$ of
infinite order has a power in $H$.
\end{Lemma}

\begin{proof}[Proof of the lemma]
If $H$ has finite index there is a partition of $G$ into a finite
number of classes modulo $H$.
The powers of $g$, in countable number, are distributed in these
classes, so there exist distinct integers $i$ and $j$ such that $g^i$
and $g^j$ are in the same class, and then $g^{j-i}\in H$.
\end{proof}

Apply this lemma to $G=\SL(2,\Z)$ and $H$ the common subgroup to $G$
and $\Gamma$, of finite index in both $G$ and $\Gamma$, and
$g=\mat{1}{0}{1}{1}$ or $g=\mat{1}{1}{0}{1}$.

We now prove Proposition~\ref{prn:cont:parab:squ:tiled}.

Since $\mat{1}{0}{m}{1}\in\Gamma$, the horizontal direction is
parabolic, so $S$ decomposes into horizontal cylinders
$C_i^\mathrm{h}$ of rational moduli.
Replacing $\mat{1}{0}{m}{1}$ with one of its powers if necessary,
suppose it fixes the boundaries of these cylinders.
This means their moduli are multiples of $1/m$.
Calling $w_i^\mathrm{h}$, $h_i^\mathrm{h}$ the widths and heights of
these cylinders, we have relations $h_i^\mathrm{h} / w_i^\mathrm{h} =
k_i / m$ for some integers $k_i$.

By a similar argument, since $\mat{1}{n}{0}{1}\in\Gamma$, the vertical
direction is also parabolic, and $S$ decomposes into vertical
cylinders $C_j^\mathrm{v}$ of rational moduli $h_j^\mathrm{v} /
w_j^\mathrm{v}= k'_j / n$ for some integers $k'_j$.

Combining these two decompositions yields a decomposition of $S$ into
rectangles of dimensions $h_j^\mathrm{v} \times h_i^\mathrm{h}$ (these
rectangles are the connected components of the intersections of the
horizontal and vertical cylinders).
Here we keep on with the convention of \S\,\ref{coord:sts} about
heights and widths of cylinders.

What we want to show is that these rectangles have rational dimensions
(up to a common real scaling factor), in order to prove that $S$ is a
covering of a square torus; indeed, if the rectangles are such, then
they can be divided into equal squares, so we obtain a covering of a
square torus.
Since singular points of $S$ lie on the edges both of horizontal and
of vertical cylinders, they are at corners of rectangles and hence of
squares of the tiling, so that the covering is ramified over only one
point.

Because the cylinders in the decompositions above are made up of these
rectangles, we have $w_i^\mathrm{h} = \sum m_{ij}h_j^\mathrm{v}$ and
$w_j^\mathrm{v} = \sum n_{ji}h_i^\mathrm{h}$, where $m_{ij}, n_{ji}
\in \N$.

Combining equations, $m h_i^h = \sum k_i m_{ij} h_j^v$ and
$n h_j^v = \sum k'_j n_{ji} h_i^h$.

Then, setting $X^\mathrm{h} = (h_i^\mathrm{h})$, $X^\mathrm{v} =
(h_j^\mathrm{v})$, $M = (k_im_{ij})_{ij}$, $N = (k'_jn_{ji})_{ji}$, we
have $mX^\mathrm{h} = M X^\mathrm{v}$ and $n X^\mathrm{v} = N
X^\mathrm{h}$, so that $MN X^\mathrm{h} = mn X^\mathrm{h}$ and $NM
X^\mathrm{v} = nm X^\mathrm{v}$.

$M$, $N$ and their products are matrices with nonnegative integer
coefficients.
In view of applying the Perron--Frobenius theorem, we show that MN and
NM have powers with all coefficients positive.

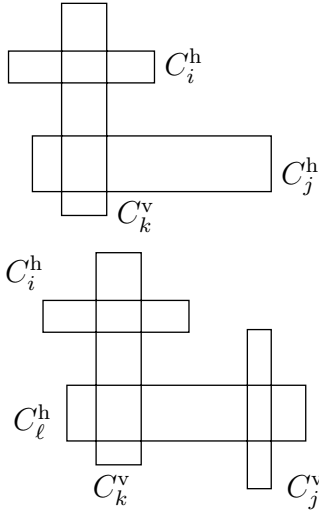
\begin{wrapfigure}{l}{128pt}
  \psfrag{Cih}{$C_i^\mathrm{h}$}
  \psfrag{Cjh}{$C_j^\mathrm{h}$}
  \psfrag{Cjv}{$C_j^\mathrm{v}$}
  \psfrag{Ckv}{$C_k^\mathrm{v}$}
  \psfrag{Clh}{$C_\ell^\mathrm{h}$}
  \includegraphics{cylguju1.eps}\\
  \includegraphics{cylguju2.eps}
\end{wrapfigure}

This results from the connectedness of $S$ and the following
observation: $M_{ij} \neq 0$ if and only if $C_i^\mathrm{h}$ and
$C_j^\mathrm{v}$ intersect; $(MN)_{ij} \neq 0$ if and only if there
exists a cylinder $C_k^\mathrm{v}$ which intersects both
$C_i^\mathrm{h}$ and $C_j^\mathrm{h}$, as in the picture; more
generally the element $i,j$ of a product of alternately $M$ and $N$
matrices is nonzero if and only if there exists a corresponding
sequence of alternately horizontal and vertical cylinders such that
two successive cylinders intersect.
So $MN$ and $NM$ do have powers with all coefficients positive.

$X^\mathrm{h}$ (resp.\ $X^\mathrm{v}$) is an eigenvector for the
eigenvalue $nm$ of the square matrix $MN$ (resp.\ $NM$).
By the Perron--Frobenius theorem, there exists a unique eigenvector
associated with the real positive eigenvalue $nm$ for the matrix $NM$
(resp.\ $MN$).
Since both matrices have rational coefficients and the eigenvalue is
rational, there exist eigenvectors with rational coefficients.
Up to scaling, they are unique by the Perron--Frobenius theorem.
This allows to conclude that $X^\mathrm{h}$ is a multiple of a vector
with rational coordinates.
From the equation $n X^\mathrm{v}=N X^\mathrm{h}$, we then conclude
that the rectangles have rational moduli and can be tiled by identical
squares.
This completes the proof of the theorem.

\subsection{A corollary}

The following result of \cite{GuHuSc} arises as a corollary of
\S\,\ref{sec:sts:arith} and
Proposition~\ref{prn:cont:parab:squ:tiled}.

\begin{Corollary}
If a subgroup $\Gamma <\SL(2,\Z)$ contains two elements
$\mat{1}{0}{m}{1}$ and $\mat{1}{n}{0}{1}$ and has infinite index in
$\SL(2,\Z)$, then $\Gamma$ cannot be realized as the Veech group of a
translation surface.
\end{Corollary}


\begin{thebibliography}{Z}
%
%
\bibitem[Ca]{Ca}
%
K.~Calta.
%
Veech surfaces and complete periodicity in genus 2.
%
Preprint. {\Small \verb+arXiv:math.DS/0205163+}
%
%
\bibitem[CuPa]{CuPa}
%
C.~J.~Cummins, S.~Pauli.
%
Congruence subgroups of PSL(2,Z) of genus less than or equal to 24.
%
\textit{Experimental mathematics} {\bf 12}:2 (2003), 243--255.\\
%
{\Small \verb+http://www.math.tu-berlin.de/~pauli/congruence/+}
%
%
\bibitem[EsMaSc]{EsMaSc}
%
A.~Eskin, H.~Masur, M.~Schmoll.
%
Billiards in rectangles with barriers.
%
\textit{Duke Math.\ J.}\ \textbf{118}:3 (2003) 427--463.
%
%
\bibitem[EsOk]{EsOk}
%
A.~Eskin, A.~Okounkov.
%
Asymptotics of numbers of branched coverings of a torus and volumes of
moduli spaces of holomorphic differentials.
%
\textit{Invent.\ Math.}\ {\bf 145}:1 (2001), 59--103.
%
%
\bibitem[Gu]{Gu}
%
E.~Gutkin.
%
Billiards on almost integrable polyhedral surfaces.
%
\textit{Ergodic Theory Dynam.\ Systems} \textbf{4}:4 (1984) 569--584.
%
%
\bibitem[GuHuSc]{GuHuSc}
%
E.~Gutkin, P.~Hubert, T.~Schmidt.
%
Affine diffeomorphisms of translation surfaces: periodic points,
fuchsian groups, and arithmeticity.
%
To appear in \textit{Ann.\ Sci.\ \'Ecole Norm.\ Sup.\ (4).}
%
%
\bibitem[GuJu1]{GuJu96}
%
E.~Gutkin, C.~Judge.
%
The geometry and arithmetic of translation surfaces with
applications to polygonal billiards.
%
\textit{Math.\ Res.\ Lett.}\ {\bf 3}:3 (1996), 391--403.
%
%
\bibitem[GuJu2]{GuJu00}
%
E.~Gutkin, C.~Judge.
%
Affine mappings of translation surfaces: geometry and arithmetic.
%
\textit{Duke Math.\ J.}\ {\bf 103}:2 (2000), 191--213.
%
%
\bibitem[Hu]{Hu}
%
P.~Hubert.
%
\textit{\'Etude g\'eom\'etrique et combinatoire de syst\`emes
dynamiques d'entropie nulle.}
%
Habilitation \`a diriger des recherches. 2002.
%
%
\bibitem[HL]{HL}
%
P.~Hubert, S.~Leli\`evre.
%
Noncongruence subgroups in $\cH(2)$.
%
To appear in \textit{Internat.\ Math.\ Res.\ Notices}.
%
%
\bibitem[HuSc00]{HuSc00}
%
P.~Hubert, T.~Schmidt.
%
Veech groups and polygonal coverings.
%
\textit{J.\ Geom.\ Physics} \textbf{35}:1 (2000), 75--91.
%
%
\bibitem[HuSc01]{HuSc01}
%
P.~Hubert, T.~Schmidt.
%
Invariants of translation surfaces.
%
\textit{Ann.\ Inst.\ Fourier (Grenoble)} {\bf 51}:2 (2001), 461--495.
%
%
\bibitem[Ka]{Ka}
%
E.~Kani.
%
The number of genus 2 covers of an elliptic curve.
%
Preprint (2003).
%
%
\bibitem[KeMaSm]{KeMaSm}
%
S.~Kerckhoff, H.~Masur, J.~Smillie.
%
Ergodicity of billiard flows and quadratic differentials.
%
\textit{Ann.\ of Math.\ (2)} \textbf{124}:2 (1986), 293--311.
%
%
\bibitem[KeSm]{KeSm}
%
R.~Kenyon, J.~Smillie.
%
Billiards on rational-angled triangles.
%
\textit{Comment.\ Math.\ Helv.}\ {\bf 75}:1 (2000), 65--108.
%
%
\bibitem[KoZo]{KoZo}
%
M.~Kontsevich, A.~Zorich.
%
Connected components of the moduli space of holomorphic differentials
with prescribed singularities.
%
\textit{Invent.\ Math.}\  \textbf{153}:3 (2003), 631--678.
%
%
\bibitem[Ma]{Ma}
%
H.~Masur.
%
Interval exchange transformations and measured foliations.
%
\textit{Ann.\ of Math.\ (2)} \textbf{115}:1 (1982), 169--200.
%
%
\bibitem[Mc]{Mc}
%
C.~T.~McMullen.
%
Billiards and Teich\-m\"uller curves on Hilbert modular surfaces.
%
\textit{J.\ Amer.\ Math.\ Soc.}\ \textbf{16}:4 (2003), 857--885.
%
%
\bibitem[Mc2]{Mc2}
%
C.~T.~McMullen.
%
Teich\-m\"uller curves in genus two: discriminant and spin.
%
Preprint (2004).
%
%
\bibitem[M\"o]{Moe}
%
M.~M\"oller.
%
Teich\-m\"uller curves, Galois actions and $\widehat{GT}$-relations.
%
Preprint (2003). {\Small \verb+arXiv:math.AG/0311308+}
%
%
\bibitem[Ri]{Ri}
%
J.~Rivat.
%
Private communication.
%
%
\bibitem[Schmi]{Schmi}
%
G.~Schmith\"usen.
%
An algorithm for finding the Veech group of an origami.
%
To appear in \textit{Experimental Mathematics}.
%
%
\bibitem[Schmo]{Schmo}
%
M.~Schmoll.
%
On the asymptotic quadratic growth rate of saddle connections and
periodic orbits on marked flat tori.
%
\textit{Geom.\ Funct.\ Anal.}\ \textbf{12}:3 (2002), 622--649.
%
%
\bibitem[Th]{Th}
%
W.~P.~Thurston.
%
On the geometry and dynamics of diffeomorphisms of surfaces.
%
\textit{Bull.\ Amer.\ Math.\ Soc.\ (N.S.)} {\bf 19}:2 (1988), 417-431.
%
%
\bibitem[Ve82]{Ve82}
%
W.~A.~Veech.
%
Gauss measures for transformations on the space of interval exchange
maps.
%
\textit{Ann.\ of Math.\ (2)} \textbf{115}:1 (1982), 201--242.
%
%
\bibitem[Ve87]{Ve87}
W.~A.~Veech.
%
Boshernitzan's criterion for unique ergodicity of an interval
exchange transformation.
%
\textit{Ergodic Theory Dynam.\ Systems} \textbf{7}:1 (1987),
149--153.
%
%
\bibitem[Ve89]{Ve89}
%
W.~A.~Veech.
%
Teich\-m\"uller curves in moduli space, Eisenstein series and an
application to triangular billiards.
%
\textit{Invent.\ Math.}\ {\bf 97}:3 (1989), 553-583.
%
%
\bibitem[Ve92]{Ve92}
%
W.~A.~Veech.
%
The billiard in a regular polygon.
%
\textit{Geom.\ Funct.\ Anal.}\ \textbf{2}:3 (1992), 341--379.
%
%
\bibitem[Ve95]{Ve95}
%
W.~A.~Veech.
Geometric realizations of hyperelliptic curves.
%
\textit{Algorithms, fractals, and dynamics (Okayama/Kyoto, 1992)},
217--226, Plenum, New York, 1995.
%
%
\bibitem[Vo]{Vo}
%
Ya.~B.~Vorobets.
%
Planar structures and billiards in rational polygons: the Veech
alternative.
%
\textit{Russ.\ Math.\ Surv.}\ {\bf 51}:5 (1996), 779-817.
%
%
\bibitem[Wa]{Wa}
%
C.~C.~Ward.
%
Calculation of fuchsian groups associated to billiards in a rational
triangle.
%
\textit{Ergodic Theory Dynam.\ Systems} \textbf{18}:4 (1998),
1019--1042.
%
%
\bibitem[Zo]{Zo}
%
A.~Zorich.
%
Square tiled surfaces and Teich\-m\"uller volumes of the moduli
spaces of abelian differentials.
%
\textit{Rigidity in dynamics and geometry (Cambridge, 2000)},
459--471, Springer, Berlin, 2002.
%
%
\end{thebibliography}
\end{document}